\documentclass[11pt]{article}

\usepackage{titlesec}
\titleformat*{\section}{\LARGE\bfseries}
\titleformat*{\subsection}{\Large\bfseries}
\titleformat*{\subsubsection}{\large\bfseries}

\makeatletter

\@addtoreset{equation}{section}

\@addtoreset{figure}{section}

\@addtoreset{table}{section}

\makeatother

\usepackage{amsmath, amsthm,amsfonts, amssymb, bm, fancybox, multirow, hhline, psfrag, mathtools, color,algorithmic,algorithm}
\usepackage[dvipdfmx]{graphicx}
\usepackage{here}
\usepackage{enumerate}
\usepackage{ascmac}
\usepackage{url}
\usepackage{color}
\theoremstyle{definition}
\newtheorem{theorem}{Theorem}[section]

\newtheorem{lemma}{Lemma}[section]

\newtheorem{assumption}{Assumption}[section]

\newcommand{\argmin}{\mathop{\rm arg~min}\limits}
\newcommand{\argmax}{\mathop{\rm arg~max}\limits}

\makeatletter
\renewcommand{\ALG@name}{アルゴリズム}
\makeatother
\allowdisplaybreaks[1]

\title{\huge{Fast same-step forecast in SUTSE model and its theoretical properties}\\[0.5cm]}

\author{      
	\LARGE{Wataru Yoshida and Kei Hirose}\\[0.5cm]
	\LARGE{Kyushu University}\\[0.5cm]
}
\date{
}

\usepackage{natbib}
\bibpunct[:]{[}{]}{,}{a}{}{,}

\begin{document}

\maketitle
\thispagestyle{empty}

\setcounter{page}{1}
\pagestyle{plain}

%%%%%%%%%%%%%%%%%%%%%%%%%%%%%%%%%%%%%%%%%%%%%%%%%%%%%%%%%%%%%%%%%%%%
\begin{abstract}
We consider the problem of forecasting multivariate time series by a Seemingly Unrelated Time Series Equations (SUTSE) model. The SUTSE model usually assumes that error variables are correlated. A crucial issue is that the model estimation requires heavy computational loads because of a large matrix computation, especially for high-dimensional data. To alleviate the computational issue, we propose a two-stage procedure for forecasting. First, we perform the Kalman filter as if error variables are uncorrelated; that is, univariate time-series analyses are conducted separately to avoid a large matrix computation. Next, the forecast value is computed by using a distribution of forecast error. The proposed algorithm is much faster than the ordinary SUTSE model because we do not require a large matrix computation. Some theoretical properties of our proposed estimator are presented. Monte Carlo simulation is performed to investigate the effectiveness of our proposed method. The usefulness of our proposed procedure is illustrated through a bus congestion data application.
\end{abstract}
\section{Introduction}
Multivariate time series data analysis has been recently developed in various fields to achieve high-quality forecasts and investigate the correlations among time series: for example, forecast of energy consumption and crowd-flow \citep{gong2020online}. In particular, a state-space model is applied widely for forecasting time series, and a number of model estimation procedures have been proposed. One of the most famous methods is the Kalman filter \citep{kalman1960new} and its extensions: for example, the Adaptive Kalman filter \citep{mohamed1999adaptive} and the Robust Kalman filter \citep{koch1998robust}. In recent times, the deep Kalman filtering network which fuses deep neural networks and the Kalman filter, has been proposed as well \citep{lu2018deep}. The Kalman filter and its extensions have been used in various fields of research, such as tracking \citep{jondhale2018kalman}, photovoltaic forecasting \citep{pelland2013solar}, and traffic volume forecasting \citep{xie2007short}.

In this study, we employ a Seemingly Unrelated Time Series Equations (SUTSE) model \citep{doi:10.1080/07350015.1990.10509778,antoniou2013state}, a special case of the state-space model. In the SUTSE model, multiple univariate time series equations are combined to express a single multivariate linear Gaussian state-space model. The SUTSE model usually assumes that components of a noise vector are correlated. In this case, the components of an observation vector are also correlated. As an example of an analysis using such a correlation structure, \citet{moauro2005temporal} employed temporal disaggregation; that is, a low-frequency time series is transformed into a high-frequency series by interpolation. They performed the interpolation by making use of the correlation structure of a multivariate time series with different frequencies.

The SUTSE model can be applied not only to interpolation but also to forecasting. Suppose we have $n$ observations of $d$-dimensional time series data, $\{\{y_1,\dots, y_n\}: y_t=(y_{1,t},\dots,y_{d,t})^T,  ~ t = 1,\dots , n \}$, and some components of $y_{n+1}$, say $y_{\mathcal{A}, n+1}$. Here, $\mathcal{A} \subset \{1,\dots, d\}$ and $y_{\mathcal{A} ,n+1}$ indicates a subvector of $y_{n+1}$ whose indices consist of $\mathcal{A}$. We may consider two types of forecasts: the one-step-ahead forecast and the same-step forecast. The one-step-ahead forecast calculates the forecast value $y_{n+1}$ by using the ${y_1,...,y_n}$. In the same-step forecast, we forecast the value of $y_{k,n+1}$ by using $\{y_1,\dots, y_n\}$ and $y_{\mathcal{A}, n+1}$, where $k \notin \mathcal{A}$.

\begin{figure}[H]
\centering
\includegraphics[width=12cm, bb=0 0 532 292]{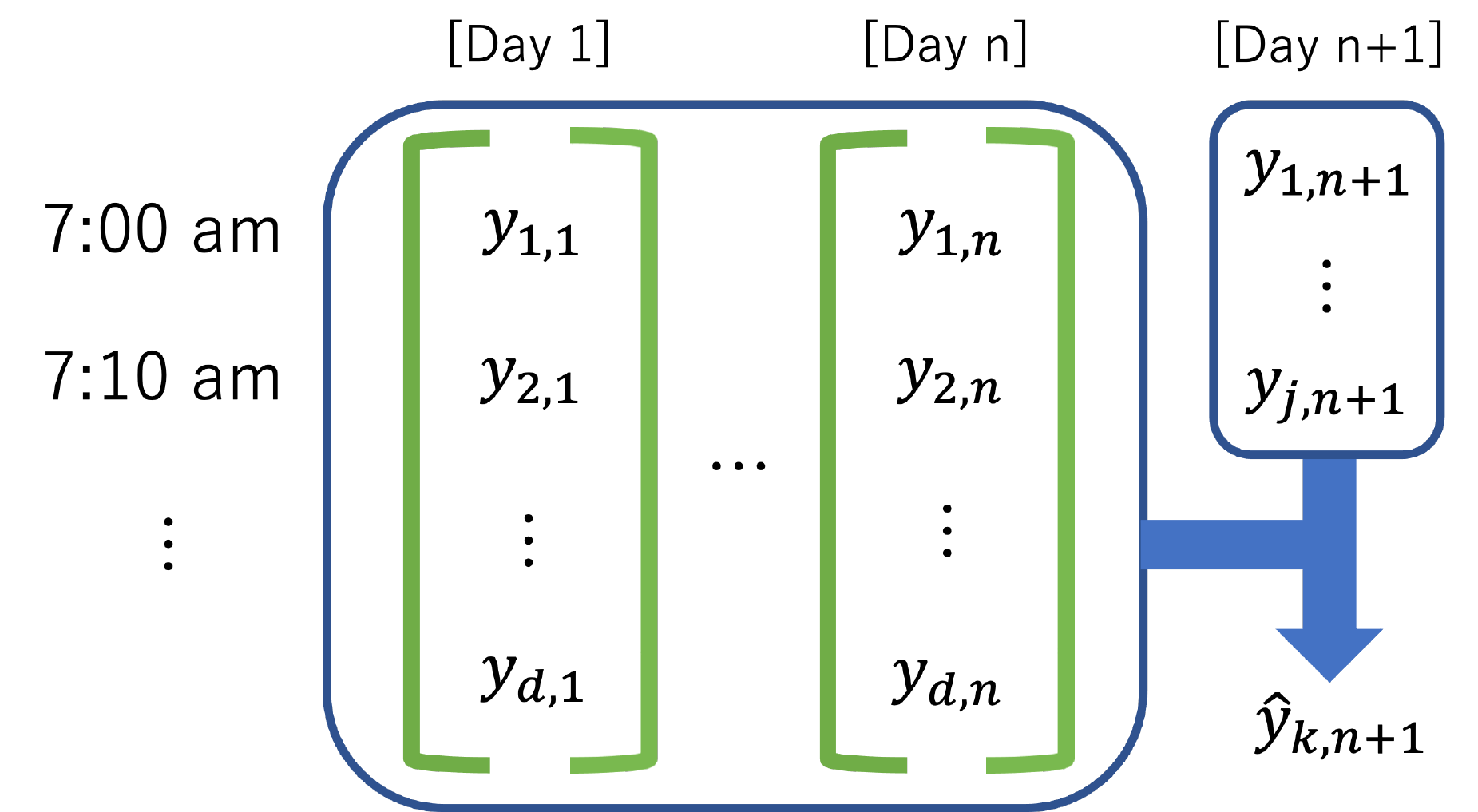}
\caption{Example: bus congestion same-step forecast}
\label{fig:zu13}
\end{figure}

Figure \ref{fig:zu13} shows one example of the application of the same-step forecast. In this example, $y_{i,t}$ indicates the congestion of $i$th bus on $t$th day. We forecast the congestion using past congestion data; in other words, $\mathcal{A} = \{1,\dots,j \}$ when we forecast $y_{k,n+1}$, where $k$ satisfies $k > j$. The same-step forecast can be used not only for a bus congestion forecast but also for a wide variety of practical applications, including electricity demand and price forecasting.

The same-step forecast is performed by modifying the one-step forecast based on the covariance matrix estimation of the one-step-ahead forecast error. This estimation is conducted with the Kalman Filter. However, the Kalman filter requires heavy computational loads because of a large matrix computation, especially for high-dimensional data. The computational complexity is about $O(\text{max}\{ d^{3}, p^{3} \} )$ \citep{willner1976kalman}, where $p$ denotes the dimension of the state vector. In fact, we attempted to analyze 32-dimensional bus congestion data with the SUTSE model. Despite using the package FKF \citep{FKF}, which can perform the Kalman filter very fast because of the implementation in C \citep{tusell2011kalman}, it took approximately 11 hours to complete the analysis. The details of this experiment are given in Section \ref{sec:bus}.

One way to handle this computational issue would be to apply a method that accelerates the Kalman Filter algorithm. For a general linear Gaussian state-space model, a method to accelerate the Kalman filter is proposed by \citet{koopman2000fast}. This method transforms the observation vectors into a univariate time series and then applies the Kalman filter to this univariate time series. A transformed univariate time series results in larger sample sizes than the original observation vectors, and its state vectors have the same dimension as the original vectors. Thus, this method reduces the cost of matrix calculations related to the original observation vectors. However, the matrix calculations related to the state vectors of transformed univariate time series must be done more times than the ordinary Kalman filter. Empirically, this method works well when $p\leq d$. Meanwhile, in the SUTSE model, it is mostly $p>d$. Therefore, this method may not speed up the Kalman filter adequately.

As seen above, the general methods for accelerating the Kalman filter in state-space models may not always speed up the SUTSE model. Thus, in this study, we propose a simple and faster method specialized for the same-step forecast in the SUTSE model. Specifically, the following two-stage procedure is proposed. First, a model estimation with the Kalman filter is performed separately for each dimension of observation vectors, as if components of the observation vector are uncorrelated, and the one-step-ahead forecast value is computed. With this procedure, the cost of matrix calculations involving both the observation and state vectors in the Kalman filter is significantly reduced. In addition, parallel computing can be applied. Next, the mean vector and the covariance matrix of the one-step-ahead forecast error are estimated using the results obtained from the Kalman filter in the first step. The same-step forecast value is computed by using the one-step-ahead forecast error distribution. Note that the mean vector and covariance matrix of the one-step-ahead forecast error are time-varying, and thus they are viewed as time-varying parameters. We show that these time-varying parameters converge under some assumptions. In particular, the mean vector converges to 0; consequently, we estimate the mean vector as 0. The sample covariance matrix is shown to be the consistent estimator of the limiting value; hence, we use the sample covariance matrix as an estimator of the covariance matrix of the one-step-ahead forecast error. Another possible covariance estimation method is applying the graphical lasso \citep{friedman2008sparse} to the sample covariance matrix. When the sample size is not large enough compared with $d$, this estimation would be more suitable. Regardless of which covariance estimation method is chosen, the second step does not take long. A Monte Carlo simulation shows that the proposed algorithm slightly sacrifices the forecast accuracy, while the computational time is greatly improved, resulting in a practical method of the same-step forecast for high-dimensional time series data. We also apply our proposed method to the bus congestion data. The result shows that the proposed method took about 8 seconds to analyze, while the existing method took approximately 11 hours.

The rest of this paper is structured as follows. In Section 2, the forecast methods using the SUTSE model are presented. We also review the Kalman filter in this section. In Section 3, we propose the fast same-step forecasting method and provide its theoretical properties. In Section 4, we prove the convergence of the mean vector and the covariance matrix of the one-step-ahead forecast error. Using this convergence, in Section 5, we prove the consistency of the estimator in the second step of the fast method. In Section 6, we conduct a Monte Carlo simulation and bus congestion forecasting to investigate the computational time and the forecast accuracy of our proposed algorithm. Section 7 presents the conclusion, and the Appendix contains proofs of lemmas and other supplemental work.
%%%%%%%%%%%%%%%%%%%%%%%%%%%%%%%%%%%%%%%%%%%%%%%%%%%%%%%%%%%%%%%%%%%%
\section{Forecasts using the SUTSE model}
\subsection{SUTSE model}
Let $y_1,y_2,\dots,y_n$ be $d\times 1$ observation vectors with $y_{t}:=(y_{1,t},\dots ,y_{d,t})^{T}$. Assume that $y_{j,t}$ follow the linear-Gaussian state-space model:
\begin{align} \label{eq:SUTSE1}
\left\{ \begin{array}{ll}
 y_{j,t} = Z_{t}^{(j)}\alpha _{t}^{(j)}+\varepsilon _{j,t}   & \\
 \alpha _{t+1}^{(j)} = T_{t}^{(j)}\alpha _{t}^{(j)} + \eta _{t}^{(j)}  & \end{array} \right.  ( t = 1,2,\dots,n, ~ j = 1,2,\dots, d),
\end{align}
where $\alpha _{t}^{(j)}$ is a $p^{(j)}\times 1$ state vector, $\varepsilon _{j,t}$ is an observation noise, $\eta _{t}^{(j)}$ is a $p^{(j)}\times 1$ state noise vector, $Z_{t}^{(j)}$ is a $1\times p^{(j)}$ design matrix, and $T_{t}^{(j)}$ is $p^{(j)}\times p^{(j)}$ a transition matrix. Let $Z_{t}:=Diag(Z_{t}^{(1)},\dots ,Z_{t}^{(d)})$, $\alpha_{t}:=(\alpha _{t}^{(1)T},\dots ,\alpha_{t}^{(d)T})^T$, $T_{t}:=Diag(T_{t}^{(1)},\dots , T_{t}^{(d)})$, $\eta_{t}:=(\eta_{t}^{(1)T},\dots ,\eta_{t}^{(d)T})^T$, and $p:=p^{(1)}+\dots p^{(d)}$. Then model \eqref{eq:SUTSE1} can be rewritten as the form of the multivariate linear-Gaussian state-space model:
\begin{align} \label{eq:model}
\left\{ \begin{array}{ll}
 y_{t} = Z_{t}\alpha _{t}+\varepsilon _{t}   & \\
 \alpha _{t+1} = T_{t}\alpha _{t}+\eta _{t}  & \end{array} \right.  ( t = 1,2,\dots,n),
\end{align}
where $\alpha _{t}$ is a $p\times 1$ state vector, $\varepsilon _{t}$ is a $d\times 1$ observation noise vector, $\eta _{t}$ is a $p\times 1$ state noise vector, $Z_{t}$ is a $d\times p$ design matrix, and $T_{t}$ is a $p\times p$ transition matrix. We assume $\varepsilon _{t}\sim N\left( 0,\Sigma _{\varepsilon }\right)$, $\eta_{t}\sim N\left( 0,\Sigma_{\eta} \right)$, $\alpha _{1}\sim N\left( a_{1},P_{1}\right)$, and $\varepsilon _{t}$, $\varepsilon _{k}(t\neq k)$, $\eta_{t}$, and $\eta_{k}$ are mutually uncorrelated. This model is called Seemingly Unrelated Time Series Equations model (SUTSE model) \citep{doi:10.1080/07350015.1990.10509778}.

\subsection{One-step-ahead forecast}\label{sec:tommorow_forecast}
We now consider the forecast of $y_{n+1}$ when $y_1,y_2,\dots,y_n$ are given. The Kalman filter is well known to be used for this forecast. Let $a_{t}:=E(\alpha_{t}|Y_{t-1})$, $P_{t}:=V(\alpha_{t}|Y_{t-1})$, $v_{t}:=y_{t}-Z_{t}a_{t}$, and $F_{t}:=V(v_{t}|Y_{t-1})$, where $Y_{t}$ denotes $\{ y_1,y_2,\dots,y_t  \}$. These conditional expectations and covariance matrices can be computed with the Kalman filter as follows:
\begin{align} \label{eq:kalman}
\left\{ \begin{array}{ll}
v_{t} = y_{t}-Z_{t}a_{t} 
\\
F_{t} = Z_{t}P_{t}Z_{t}^{\mathrm{T}}+\Sigma _{\varepsilon }
\\
a_{t+1} = T_{t}a_{t}+T_{t}K_{t}v_{t}
\\
P_{t+1} = T_{t}P_{t}L_{t}^{\mathrm{T}}T_{t}^{\mathrm{T}}+\Sigma _{\eta}  \end{array} \right.  ( t = 1,2,\dots,n),
\end{align}
where $K_{t}:=P_{t}Z_{t}^{\mathrm{T}}F_{t}^{\mathrm{-1}}$ and $L_{t}:=I_{p}-K_{t}Z_{t}$. Note that $a_{1}$ and $P_{1}$ are the mean vector and the covariance matrix of the initial state vector $\alpha _{1}$ respectively. Then, we define ``one-step-ahead forecast'' of $y_{t}$ as
\begin{align}\label{eq:tommorow_forecast}
\bar{y}_{t} &:= Z_{t}a_{t},
\end{align}
and call $v_{t}$ ``one-step-ahead forecast error'' in this paper. In particular, we can get the one-step-ahead forecast of $y_{n+1}$ as $\bar{y}_{n+1} = Z_{n+1}a_{n+1}$.
Note that by the definition of $a_{t}$, we have $\bar{y}_{t}=E(Z_{t}\alpha_{t}|Y_{t-1})=E(y_{t}|Y_{t-1})$ and $E(v_{t}|Y_{t-1})=E(y_{t}-Z_{t}a_{t}|Y_{t-1})=E(y_{t}-\bar{y}_{t}|Y_{t-1})=0$.

\subsection{Same-step forecast}\label{sec:today_forecast}
We now consider the forecast of $y_{k,n+1}$ when $y_1,y_2,\dots,y_n$ and also $y_{1,n+1},\dots, y_{j,n+1}$ are given with $k>j$. The one-step-ahead forecast of $y_{k,n+1}$ can be modified by using the conditional expectation of the one-step-ahead forecast error $v_{k,n+1}$ as follows:
\begin{align}\label{eq:today_forecast}
\hat{y}_{k,n+1} &= [Z_{n+1}a_{n+1}](k)+E\Bigl( v_{k,n+1}|Y_{n},v_{n+1}(1:j) \Bigl) ,
\end{align}
where $[Z_{n+1}a_{n+1}](k)$ denotes the $k$th component of $Z_{n+1}a_{n+1}$ and $c(i:j)$ denotes $(c_{i},\dots , c_{j})^{T}$. We call this forecast ``same-step forecast''. To compute the modification term (the second term) of the right-hand side of \eqref{eq:today_forecast}, we use the formula for a conditional expectation of a multivariate normal distribution (e.g., see \citep{eaton1983multivariate}):
\begin{align}\label{fomula_muln}
&E(x_{1}|x_{2})=E(x_{1})+{\rm Cov}(x_{1},x_{2})V(x_{2})^{-1}(x_{2}-E(x_{2})), \nonumber
\\
&\text{where $x_{1}$ and $x_{2}$ follow multivariate normal distribution.}
\end{align}
Hence, the modification term can be expressed as
\begin{align*}
&E\Bigl( v_{k,n+1}|Y_{n},v_{n+1}(1:j) \Bigl)
\\
&=E(v_{k,n+1}|Y_{n}) + {\rm Cov} ( v_{k,n+1},v_{n+1}(1:j)|Y_{n} ) V ( v_{n+1}(1:j)|Y_{n} )^{-1} \Bigl\{ v_{n+1}(1:j) - E ( v_{n+1}(1:j) | Y_{n} ) \Bigl\} .
\end{align*}
Since $E(v_{n+1}|Y_{n})=0$ and $V(v_{n+1}|Y_{n})=F_{n+1}$, we have
\begin{align}\label{eq:syuuseikou}
E\Bigl( v_{k,n+1}|Y_{n}, v_{n+1}(1:j) \Bigl) =F_{n+1}(k,1:j)F_{n+1}(1:j,1:j)^{-1} v_{n+1}(1:j) .
\end{align}
Here, $A(k_{1}:k_{2},l_{1}:l_{2})$ denotes the submatrix consisting of rows of indices from $k_{1}$ to $k_{2}$ and columns of indices from $l_{1}$ to $l_{2}$ of a matrix $A$. Substituting \eqref{eq:syuuseikou} into \eqref{eq:today_forecast}, the same-step forecast can be computed as
\begin{eqnarray}\label{eq:today_forecast_cal}
\hat{y}_{k,n+1}= [Z_{n+1}a_{n+1}](k)+F_{n+1}(k,1:j)F_{n+1}(1:j,1:j)^{-1} v_{n+1}(1:j) .
\end{eqnarray}
Now, the modification term $E\Bigl( v_{k,n+1}|Y_{n}, v_{n+1}(1:j) \Bigl)$ satisfies
\begin{eqnarray}
E\Bigl( v_{k,n+1}|Y_{n}, v_{n+1}(1:j) \Bigl) = E\Bigl( v_{k,n+1}|v_{n+1}(1:j) \Bigl) .
\end{eqnarray}
The above equation is proved from the following facts: 
\begin{align*}
&E(v_{n+1})=E\Bigl( E(v_{n+1} | Y_{n})\Bigl) =0,
\\
&V(v_{n+1})=E\Bigl( V(v_{n+1}|Y_{n})\Bigl) + V\Bigl( E(v_{n+1}|Y_{n})\Bigl) =E(F_{n+1})+V(0)=F_{n+1}.
\end{align*}
It follows that
\begin{align}
&E\Bigl( v_{k,n+1}| v_{n+1}(1:j) \Bigl) \nonumber
\\
&=E(v_{k,n+1})+{\rm Cov}\Bigl( v_{k,n+1}, v_{n+1}(1:j) \Bigl) V\Bigl( v_{n+1}(1:j) \Bigl) ^{-1} \Bigl\{ v_{n+1}(1:j) -E\Bigl( v_{n+1}(1:j) \Bigl) \Bigl\} \nonumber
\\
&=F_{n+1}(j+1,1:j)F_{n+1}(1:j,1:j)^{-1} v_{n+1}(1:j) = E\Bigl( v_{k,n+1}|Y_{n},v_{n+1}(1:j) \Bigl), \nonumber
\end{align}
where we used \eqref{fomula_muln} in the first equation. Therefore, the same-step forecast \eqref{eq:today_forecast} can be rewritten as follows:
\begin{align}\label{eq:today_forecast_1}
\hat{y}_{k,n+1} &= [Z_{n+1}a_{n+1}](k)+E\Bigl( v_{k,n+1}|v_{n+1}(1:j) \Bigl).
\end{align}
Clearly, the right hand side of \eqref{eq:today_forecast_1} is equal to that of \eqref{eq:today_forecast_cal}.
%%%%%%%%%%%%%%%%%%%%%%%%%%%%%%%%%%%%%%%%%%%%%%%%%%%%%%%%%%%%%%%%%%%%
\section{Fast method}
In practice, the estimation of unknown parameters is required. In the SUTSE model, the maximum likelihood estimation is usually performed:
\begin{align}\label{yuudo_existing}
& \hat{\Theta } = \argmax_{\Theta } \sum^n_{t=1} \left\{ -\frac{1}{2}v_{t}^T F_{t}^{-1} v_{t} - \frac{1}{2} \log (|F_{t}|) \right\},
\end{align}
where $\Theta$ denotes unknown parameters. The maximum likelihood estimator is not explicitly obtained; thus, numerical optimization, such as the quasi-Newton method, is used. Here, $F_{t}$ and $v_{t}$ in \eqref{yuudo_existing} are obtained by the Kalman filter \eqref{eq:kalman}. Therefore, in parameter estimation, the Kalman filter must be run many times until numerical optimization converges. Especially in high dimensions (i.e., $d$ or $p$ is large), iterative calculations of large matrix multiplications in the Kalman filter make parameter estimation computationally expensive.
\subsection{Kalman filter with disregarding correlation}\label{sec:kalman_nocor}
To address the problem mentioned above, we conduct a univariate time series analysis including a separate parameter estimation. In this method, the Kalman filter is performed separately for each dimension as follows:
\begin{align} \label{eq:kalman_nocor1}
\left\{ \begin{array}{ll}
v_{j,t}' = y_{j,t}-Z_{t}^{(j)}a_{t}'^{(j)} 
\\
F_{t}'^{(j)} = Z_{t}^{(j)}P_{t}'^{(j)}Z_{t}^{(j)T}+\Sigma _{\varepsilon }(j,j)
\\
a_{t+1}'^{(j)} = T_{t}^{(j)}a_{t}'^{(j)}+T_{t}^{(j)}K_{t}'^{(j)}v_{t}'^{(j)}
\\
P_{t+1}'^{(j)} = T_{t}^{(j)}P_{t}'^{(j)}L_{t}'^{(j)T}T_{t}^{(j)T}+\Sigma _{\eta}^{(j)}
\end{array} \right.  ( t = 1,2,\dots,n, ~ j = 1,2,\dots, d)
\end{align}
where $K_{t}'^{(j)}:=P_{t}'^{(j)}Z_{t}^{(j)T}F_{t}'^{(j)-1}$, $L_{t}'^{(j)}:=I_{p^{(j)}}-K_{t}'^{(j)}Z_{t}^{(j)}$, and $A(i,j)$ denotes $(i,j)$th component of $A$.
This method reduces the cost of matrix calculations compared with \eqref{eq:kalman}. For example, the cost of deriving $F_{t}^{-1}$ in \eqref{eq:kalman} is approximately $O(d^{3})$, but the computations for $F_{t}'^{(1)-1},\dots ,F_{t}'^{(d)-1}$ in \eqref{eq:kalman_nocor1} require only $O(d)$. In practice, parameter estimation is performed for each dimension as follows:
\begin{align}\label{yuudo_fast}
& \hat{\Theta }'^{(j)} = \argmax_{\Theta'^{(j)} } \sum^n_{t=1} \left\{ -\frac{v_{j,t}'^{2}}{2 F_{t}'^{(j)} }  - \frac{1}{2} \log ( F_{t}'^{(j)} ) \right\},
\end{align}
where $\Theta'^{(j)}$ denotes unknown parameters related to $y_{j,t}$. This method also reduces time-consuming estimation. Only the estimation of $\Sigma _{\varepsilon }(j,j)$ and $\Sigma _{\eta}^{(j)}$ for each $j=1,\dots,d$ is required instead of $\Sigma _{\varepsilon }$ and $\Sigma _{\eta}$. In addition, the parameter estimation \eqref{yuudo_fast} can be performed in parallel for each $j = 1,2,\dots, d$.

Let $v_{t}':=(v_{1,t}',\dots ,v_{d,t}')^{T}$, $a_{t}':=(a_{t}'^{(1)},\dots ,a_{t}'^{(d)})^{T}$, $F_{t}':=Diag(F_{t}'^{(1)},\dots ,F_{t}'^{(d)})$, $P_{t}':=Diag(P_{t}'^{(1)},\dots ,P_{t}'^{(d)})$, $K_{t}':=Diag(K_{t}'^{(1)},\dots ,K_{t}'^{(d)})$, $L_{t}':=Diag(L_{t}'^{(1)},\dots ,L_{t}'^{(d)})$, $\Sigma _{\varepsilon }':=Diag(\Sigma _{\varepsilon }(1,1) , \dots ,\Sigma _{\varepsilon }(d,d))$, and $\Sigma _{\eta }':=Diag(\Sigma _{\eta }^{(1)} , \dots ,\Sigma _{\eta }^{(d)})$; then, \eqref{eq:kalman_nocor1} can be rewritten as follows:
\begin{align} \label{eq:kalman_nocor2}
\left\{ \begin{array}{ll}
v_{t}' = y_{t}-Z_{t}a_{t}'
\\
F_{t}' = Z_{t}P_{t}'Z_{t}^{\mathrm{T}}+\Sigma _{\varepsilon }'
\\
a_{t+1}' = T_{t}a_{t}'+T_{t}K_{t}'v_{t}'
\\
P_{t+1}' = T_{t}P_{t}'L_{t}'^{\mathrm{T}}T_{t}^{\mathrm{T}}+\Sigma _{\eta}'
\end{array} \right.  ( t = 1,2,\dots,n).
\end{align}
See \eqref{eq:kalman_nocor2}. If $P_{1}'=P_{1}$, $\Sigma _{\varepsilon }'=\Sigma _{\varepsilon }$, and $\Sigma _{\eta }'=\Sigma _{\eta }$, the values derived by \eqref{eq:kalman_nocor1} are equal to the values derived by \eqref{eq:kalman}. Here, for example, the equation $\Sigma _{\varepsilon }'=\Sigma _{\varepsilon }$ means that the components of the observation noise $\varepsilon _{t}$ are uncorrelated. Nevertheless, it should be assumed that such correlations actually exist (i.e., $P_{1}\neq P_{1}'$, $\Sigma _{\varepsilon }\neq \Sigma _{\varepsilon }'$, or $\Sigma _{\eta }\neq \Sigma _{\eta }'$) to consider the same-step forecast. Specifically, if such correlations do not exist, one-step-ahead forecast errors $v_{1,t},\dots ,v_{d,t}$ are mutually uncorrelated. Hence the modification term of the same-step forecast in \eqref{eq:today_forecast_1} becomes 0, and so the one-step-ahead forecast and the same-step forecast become equal. Under $P_{1}\neq P_{1}'$, $\Sigma _{\varepsilon }\neq \Sigma _{\varepsilon }'$, or $\Sigma _{\eta }\neq \Sigma _{\eta }'$, of course, the values of \eqref{eq:kalman_nocor2} are not equal to the values derived by \eqref{eq:kalman}.

\subsection{Fast one-step-ahead forecast}
When $y_1,y_2,\dots,y_n$ are given, similarly to \eqref{eq:tommorow_forecast}, we consider the one-step-ahead forecast in the fast method as
\begin{align}\label{eq:tommorow_forecast_fast}
\bar{y}_{n+1}' &= Z_{t}a_{n+1}'.
\end{align}
Note that
\begin{align}
E(v_{n+1}'|Y_{n}) &= E(y_{n+1} - Z_{t}a_{n+1}'|Y_{n}) =Z_{t}a_{n+1} - Z_{t}a_{n+1}'. \nonumber
\end{align}
Thus, it does not always hold that $E(v_{n+1}'|Y_{n})=0$ in contrast to the existing one-step-ahead forecast \eqref{eq:tommorow_forecast}. However, we show the convergence of $E(v_{n+1}')$ as follows:
\begin{align}\label{eq:Econ}
\text{Under assumption \ref{assumption}, }E(v_{t}')\to 0~~~(t\to \infty ).
\end{align}
These assumptions and the proof of this convergence are given in Section \ref{sec:convergenceEV}.

\subsection{Fast same-step forecast}
When $y_1,y_2,\dots,y_n$ and also $y_{1,n+1},\dots, y_{j,n+1}$ are given, similarly to \eqref{eq:today_forecast_1}, we consider the same-step forecast in the fast method as
\begin{align}\label{eq:today_forecast_fast}
\hat{y}_{k,n+1}' &= [ Z_{n+1}a_{n+1}' ] (k)+E\Bigl( v_{k,n+1}'|v_{n+1}'(1:j) \Bigl) .
\end{align}
Using the formula for a conditional expectation of a multivariate normal distribution \eqref{fomula_muln}, the modification term (the second term) of \eqref{eq:today_forecast_fast} can be computed as
\begin{align*}
&E\Bigl( v_{k,n+1}'| v_{n+1}'(1:j) \Bigl)
\\
&=E(v_{k,n+1}') + {\rm Cov}\Bigl( v_{k,n+1}', v_{n+1}'(1:j) \Bigl) V\Bigl( v_{n+1}'(1:j) \Bigl) ^{-1} \Bigl\{ v_{n+1}'(1:j) - E\Bigl( v_{n+1}'(1:j) \Bigl) \Bigl\} .
\end{align*}
This equation suggests that if $E(v_{n+1}')$ and $V(v_{n+1}')$ are given, the same-step forecast \eqref{eq:today_forecast_fast} can be obtained. However, it is difficult to derive $E(v_{n+1}')$ and $V(v_{n+1}')$ unless $a_{n+1}$, $F_{n+1}$, etc. are given. Note that $a_{n+1}$, $F_{n+1}$, etc. are given in the normal Kalman filter \eqref{eq:kalman}, and thus these values cannot be obtained with the fast method. Hence, we now consider the estimations of $E(v_{n+1}')$ and $V(v_{n+1}')$.

First, as we stated above, it follows that $E(v_{n+1}')$ converges to 0. In addition, we show the convergence of $V(v_{n+1}')$ as follows:
\begin{align}\label{eq:Vcon}
\text{Under $Z_{t}=Z$, $T_{t}=T$, and assumptions \ref{assumption}, } V(v_{t}')\to ^\exists{V_{v'}}~~~(t\to \infty ).
\end{align}
These assumptions and the proof of this convergence are given in Section \ref{sec:convergenceEV}. Therefore, we have $E(v_{n+1}')\simeq 0$ and $V(v_{n+1}')\simeq V_{v'}$ for large $n$. Then, we choose $n_{0}$ that is large enough to eliminate the effect around the initial value and consider the estimation of $V_{v'}$ as follows:
\begin{eqnarray}\label{eq:vhat}
\hat V_{v'}&:=&\frac{1}{n-n_{0}+1}\sum_{t=n_{0} }^{n} v_{t}'v_{t}'^T.
\end{eqnarray}
Here, the estimator $\hat V_{v'}$ is a consistent estimator of $V_{v'}$. The proof of this consistency is given in Section \ref{sec:consistencyV}. From the above, the estimators of $E(v_{n+1}')$ and $V(v_{n+1}')$ are 0 and $\hat V_{v'}$, respectively; thus we estimate $E\Bigl( v_{k,n+1}'| v_{n+1}'(1:j) \Bigl)$ as
\begin{align*}
&\hat{E} \Bigl( v_{k,n+1}'|v_{n+1}'(1:j) \Bigl)=\hat V_{v'}(k,1:j)\hat V_{v'}(1:j,1:j)^{-1} v_{n+1}'(1:j) .
\end{align*}
Substituting this into \eqref{eq:today_forecast_fast}, the same-step forecast in the fast method we propose is derived as follows:
\begin{align}
&\tilde{y}_{k,n+1} = [ Z_{n+1}a_{n+1}' ] (k)+\hat V_{v'}(k,1:j)\hat V_{v'}(1:j,1:j)^{-1} v_{n+1}'(1:j).
\end{align}

\subsection{Estimation of $V(v_{n+1}')$ with the graphical lasso}
From \eqref{eq:Econ} and \eqref{eq:Vcon}, the expectation and the covariance matrix of $v_{t}'$ converge. Therefore, we assume that $v_{n_{0}}',\dots, v_{n}'$ follow approximately the same distribution, and consider the estimation of $V(v_{n+1}')$ using the $L_{1}$ regularization method proposed by \citet{yuan2007model}. Specifically, the estimator $\hat V_{v'}^{(glasso)}$ is derived by solving the following minimization problem:
\begin{align*}
& \hat{V}_{v'}^{(glasso)} = \hat{\Omega}^{-1}, \text{ where }  \hat{\Omega} := \argmin_{\Omega } \Bigl\{ - \log |\Omega | + tr(\Omega \hat V_{v'}) + \lambda \sum_{j=1}^{d} \sum_{k=1}^{d} | \Omega(j, k) | \Bigl\}.
\end{align*}
This minimization problem can be solved by the graphical lasso \citep{friedman2008sparse}. The graphical lasso is very fast and does not need much computation time. When $n$ is not large enough compared with $d$, this estimator $\hat V_{v'}^{(glasso)}$ would be more suitable. In fact, in the bus congestion forecasting in Section \ref{sec:bus}, the forecast accuracy is improved with the graphical lasso.
%%%%%%%%%%%%%%%%%%%%%%%%%%%%%%%%%%%%%%%%%%%%%%%%%%%%%%%%%%%%%%%%%%%%%%%%%%%%%
\section{Convergence of $E(v_{t}')$ and $V(v_{t}')$}\label{sec:convergenceEV}
This paper defines the Frobenius norm as
\begin{align*}
&\| A \|_{F} := \sqrt{\sum_{i=1}^{p}\sum_{j=1}^{q} A(i,j)^{2}} = \sqrt{ tr(AA^{T}) } ,
\end{align*}
where A is $p\times q$ matrix. Additionally, we define the convergence of matrix $A_{n} \to A(n\to \infty )$ as $\| A_{n} - A \|_{F} \rightarrow 0~(n \rightarrow \infty)$.
\\

We assume the true model is \eqref{eq:model}, considering the situation when the Kalman filter is run under the misspecified parameters as in \eqref{eq:kalman_nocor2} and the values $v_{t}'$, $F_{t}'$, $ a_{t}'$, and $P_{t}'$ are obtained. Also we assume that $Z_{t}$, $T_{t}$ are time-independent, so these can be written as $Z_{t}=Z$, $T_{t}=T$. Actually, it is not necessary that $T_{t}$ and $Z_{t}$ be time-independent to prove the convergence of $E(v_{t}')$; that is, assumptions can be mild. The mild assumption is detailed in Appendix \ref{app:conE_assum}. Then, we will prove that $E(v_{t}')$ and $V(v_{t}')$ converge under some assumptions, which are as follows:
\begin{assumption}\label{assumption}
\begin{flalign*}
(1)~&\Sigma _{\varepsilon }>0,~\Sigma _{\varepsilon }'>0.
\\
(2)~&rank(Z^{T}, (ZT)^{T}, . . . , (ZT^{p-1})^{T})^{T}=p.
\\
(3)~&\text{Decompose $\Sigma _{\eta}=RQR^{T}$, where $R$ and $Q$ are $p\times r$ matrix and $r \times r$ positive matrix.}
\\
&\text{Then }rank(R, TR, . . . , T^{p-1}R)=p.
\\
&\text{Similarly, decompose $\Sigma _{\eta}'=R'Q'R'^{T}$, where $R'$ and $Q'$ are $p\times r'$ matrix and $r' \times r'$}
\\
&\text{positive matrix. Then }rank(R', TR', . . . , T^{p-1}R')=p.
\\
(4)~& {}^\exists M>0, ~ {}^\forall k\leq l \leq t,~ \|\prod_{i=k}^{l} TL'_{t-i} \|_{F}\leq M.
\end{flalign*}
\end{assumption}
Assumption 1 implies that the observation noise vector follows a non-degenerate distribution. If assumptions 2 and 3 are satisfied, the linear system \eqref{eq:model} is considered observable and controllable, respectively. The notion of observability and controllability was suggested by \citet{KALMAN1960491}. According to \citet{gilbert1963controllability}, these are used in the study of the control theory. Assumption 4 is set by us. The product of $TL'_{t}$ as written in assumption 4 appears in computational processes of $E(v_{t}')$ and $V(v_{t}')$. This product needs to be bounded to ensure the convergence of $E(v_{t}')$ and $V(v_{t}')$. Under assumptions 1,2,3, the following lemmas hold.
\begin{lemma} \label{lem:PFKLcon}
$
\\
\text{Under assumptions 1,2,3, }\text{$P_{t}\to ^\exists P$, $P_{t}'\to ^\exists P'(t\to \infty )$.}
\\
\text{And thus }F_{t}\to ^\exists F,~F_{t}'\to ^\exists F',~K_{t}\to ^\exists K,~K_{t}'\to ^\exists K',~L_{t}\to ^\exists L,~L_{t}'\to ^\exists L'.
\\
\text{Also, }^\exists M>0,\ 0<^\exists r<1\text{ such that }\|P_{t}-P\|_{F},\ \|P'_{t}-P'\|_{F},\ \|F_{t}-F\|_{F},\ \|F'_{t}-F'\|_{F},\ \|K_{t}-K\|_{F},\ \|K'_{t}-K'\|_{F},\ \|L_{t}-L\|_{F},\ \|L'_{t}-L'\|_{F}\text{ are bounded above by $Mr^{t}$ for each.}
$
\end{lemma}

\begin{lemma} \label{lem:TLeigen}
$
\\
\text{Under assumptions 1,2,3, all eigenvalues of $TL'$ are of absolute value that is less than 1.}
\\
\text{And thus $^\exists M>0,\ 0<^\exists r<1,\  \| (TL')^{n}\|_{F} \leq Mr^{n}$ (see Appendix \ref{app:TL'geocon}).}
$
\end{lemma}
These lemmas were proved in Section 6 of \citep{chui2017kalman}. Now, the following convergences of $E(v_{t}')$ and $V(v_{t}')$ hold.
\setcounter{theorem}{0}
\begin{theorem} \label{theorem:conEV}
$
\\
\text{Under assumptions 1,2,3,4, }E(v_{t}')\to 0,~V(v_{t}')\to ^\exists{V_{v'}}(t\to \infty).
$
\end{theorem}
For the proof of this theorem, we prepare more lemmas.
\begin{lemma} \label{lem:E_a}
$
\\
E(a_{t+1}-a_{t+1}') =  TL_{t}'E(a_{t}-a_{t}').
$
\end{lemma}

\begin{lemma} \label{lem:V_a}
$
\\
E \Bigl( (a_{t+1}-a_{t+1}')(a_{t+1}-a_{t+1}')^{T} \Bigl) = T(K_{t}-K_{t}')F_{t}(K_{t}-K_{t}')^{T}T^{T}+TL_{t}'E \Bigl( (a_{t}-a_{t}')(a_{t}-a_{t}')^{T} \Bigl) L_{t}'^{T}T^{T}.
$
\end{lemma}

\begin{lemma} \label{lem:sa_TL}
$
\\
\text{Under assumptions 1,2,3, }{}^\exists M>0,\ 0<{}^\exists r<1,\ \| \prod_{i=0}^{j} TL'_{t-i} - (TL')^{j+1} \|_{F}\leq (j+1)Mr^{t},
\\
{}^\exists M >0,\ 0<{}^\exists r<1,\ \| \prod_{i=0}^{j} TL'_{t-i} \|_{F} \leq M(j+2)r^{j+1}.
$
\end{lemma}

\begin{lemma} \label{lem:con_TL}
$
\\
\text{Under assumptions 1,2,3, }\prod_{i=0}^{t-1} TL_{t-i}'\to 0(t\to \infty ).
$
\end{lemma}

\begin{lemma} \label{lem:TKFKT}
$
\\
\text{Under assumptions 1,2,3, }{}^\exists M>0,\ 0<{}^\exists r<1,\ 
\\
\| T(K_{t}-K_{t}')F_{t}(K_{t}-K_{t}')^{T}T^{T} - T(K-K')F(K-K')^{T}T^{T}\|_{F} \leq Mr^{t}.
$
\end{lemma}

The proofs of these lemmas are given in Appendix \ref{app:pf_lem}. We prove theorem \ref{theorem:conEV}, using above lemmas.
\begin{proof}
First, consider the convergence of $E(v_{t}')$. Note $E(y_{t+1}|Y_{t})=Za_{t+1}$,
\begin{align*}
&E(v_{t+1}')=E(y_{t+1}-Za_{t+1}')=E\Bigl( E(y_{t+1}-Za_{t+1}'|Y_{t})\Bigl) =E(Za_{t+1}-Za_{t+1}')
\\
&=ZE(a_{t+1}-a_{t+1}')=ZTL_{t}'E(a_{t}-a_{t}')=\dots =Z(\prod_{i=0}^{t-1} TL_{t-i}')(a_{1}-a_{1}')~~~~(\because \text{lemma \ref{lem:E_a}).}
\end{align*}
By lemma \ref{lem:con_TL}, $\prod_{i=0}^{t-1} TL_{t-i}'\to 0(t\to \infty )$, so $Z(\prod_{i=0}^{t-1} TL_{t-i}')(a_{1}-a_{1}')\to 0(t\to \infty )$. Therefore, we get $E(v_{t}')\to 0(t\to \infty )$. 
\\

Next, consider the convergence of $V(v_{t}')$. By the law of total variance,
\begin{align*}
&V(v_{t+1}')=E\Bigl( V(v_{t+1}'|Y_{t}) \Bigl) + V\Bigl( E(v_{t+1}'|Y_{t}) \Bigl)
\\
&=E \Bigl( V(v_{t+1}'|Y_{t}) \Bigl) +E \Bigl( E(v_{t+1}'|Y_{t})E(v_{t+1}'|Y_{t})^{T} \Bigl) -E \Bigl( E(v_{t+1}'|Y_{t}) \Bigl) E \Bigl( E(v_{t+1}'|Y_{t}) \Bigl) ^{T}
\\
&=E \Bigl( V(v_{t+1}'|Y_{t}) \Bigl) +E \Bigl( E(v_{t+1}'|Y_{t})E(v_{t+1}'|Y_{t})^{T} \Bigl) -E(v_{t+1}')E(v_{t+1}')^{T}.
\end{align*}
The first term can be computed as $E\Bigl( V(v_{t+1}'|Y_{t})\Bigl) =F_{t+1}$, because $V(v_{t+1}'|Y_{t})=V(y_{t+1}-Za_{t+1}'|Y_{t})=V(y_{t+1}-Za_{t+1}|Y_{t})=F_{t+1}$. The second term can be computed as follows:
\begin{align*}
&E \Bigl( E(v_{t+1}'|Y_{t})E(v_{t+1}'|Y_{t})^{T} \Bigl) =E \Bigl( E(y_{t+1}-Za_{t+1}'|Y_{t})E(y_{t+1}-Za_{t+1}'|Y_{t})^{T} \Bigl)
\\
&=E \Bigl( (Za_{t+1}-Za_{t+1}')(Za_{t+1}-Za_{t+1}')^{T} \Bigl)=ZE \Bigl( (a_{t+1}-a_{t+1}')(a_{t+1}-a_{t+1}')^{T} \Bigl) Z^{T}.
\end{align*}
From the above results, we have
\begin{align}\label{V_zenbunsan}
&V(v_{t+1}')=F_{t+1} + ZE \Bigl( (a_{t+1}-a_{t+1}')(a_{t+1}-a_{t+1}')^{T} \Bigl) Z^{T} - E(v_{t+1}')E(v_{t+1}')^{T}.
\end{align}
Now, using lemma \ref{lem:V_a} repeatedly, we have
\begin{align*}
&E \Bigl( (a_{t+1}-a_{t+1}')(a_{t+1}-a_{t+1}')^{T} \Bigl)
\\
&=T(K_{t}-K_{t}')F_{t}(K_{t}-K_{t}')^{T}T^{T}+TL_{t}'E \Bigl( (a_{t}-a_{t}')(a_{t}-a_{t}')^{T} \Bigl) L_{t}'^{T}T^{T}
\\
&=T(K_{t}-K_{t}')F_{t}(K_{t}-K_{t}')^{T}T^{T}+TL_{t}'T(K_{t-1}-K_{t-1}')F_{t-1}(K_{t-1}-K_{t-1}')^{T}T^{T}L_{t}'^{T}T^{T}
\\
&+TL_{t}'TL_{t-1}'E \Bigl( (a_{t-1}-a_{t-1}')(a_{t-1}-a_{t-1}')^{T} \Bigl) L_{t-1}'^{T}T^{T}L_{t}'^{T}T^{T}
\\
&=\dots
\\
&=T(K_{t}-K_{t}')F_{t}(K_{t}-K_{t}')^{T}T^{T}
\\
&+\sum_{j=0}^{t-2}\left\{ (\prod_{i=0}^{j} TL_{t-i}')T(K_{t-j-1}-K_{t-j-1}')F_{t-j-1}(K_{t-j-1}-K_{t-j-1}')^{T}T^{T}(\prod_{i=0}^{j} TL_{t-i}')^{T}\right\}
\\
&+(\prod_{i=0}^{t-1} TL_{t-i}')(a_{1}-a_{1}')(a_{1}-a_{1}')^{T}(\prod_{i=0}^{t-1} TL_{t-i}')^{T}.
\end{align*}
Substituting this result into \eqref{V_zenbunsan}, we get
\begin{align}\label{V_sum}
&V(v_{t+1}')=F_{t+1}+ZT(K_{t}-K_{t}')F_{t}(K_{t}-K_{t}')^{T}T^{T}Z^{T} \nonumber
\\
&+Z\sum_{j=0}^{t-2}\left\{ (\prod_{i=0}^{j} TL_{t-i}')T(K_{t-j-1}-K_{t-j-1}')F_{t-j-1}(K_{t-j-1}-K_{t-j-1}')^{T}T^{T}(\prod_{i=0}^{j} TL_{t-i}')^{T}\right\} Z^{T} \nonumber
\\
&+Z(\prod_{i=0}^{t-1} TL_{t-i}')(a_{1}-a_{1}')(a_{1}-a_{1}')^{T}(\prod_{i=0}^{t-1} TL_{t-i}')^{T}Z^{T} - E(v_{t+1}')E(v_{t+1}')^{T}.
\end{align}
By lemma \ref{lem:PFKLcon}, the first and second terms of \eqref{V_sum} converge to
\begin{align*}
&F_{t+1}\to F(t\to \infty ),
\\
&ZT(K_{t}-K_{t}')F_{t}(K_{t}-K_{t}')^{T}T^{T}Z^{T}\to ZT(K-K')F(K-K')^{T}T^{T}Z^{T}(t\to \infty ).
\end{align*}
The fourth and fifth terms converge to 0, because $\prod_{i=0}^{t-1} TL_{t-i}'\to 0(t\to \infty )$ and $E(v_{t}')\to 0(t\to \infty )$ from $\prod_{i=0}^{t-1} TL_{t-i}'\to 0(t\to \infty )$(lemma \ref{lem:con_TL}). Finally, we consider the third term. We define $A_{t} := T(K_{t}-K_{t}')F_{t}(K_{t}-K_{t}')^{T}T^{T}$, $A := T(K-K')F(K-K')^{T}T^{T}$ and
\begin{align}\label{S_t}
S_{t} := \sum_{j=0}^{t-2} \left\{ (\prod_{i=0}^{j} TL_{t-i}')A_{t-j-1}(\prod_{i=0}^{j} TL_{t-i}')^{T}\right\},\ S^{*}_{t} := \sum_{j=0}^{t-2} \left\{ (TL')^{j+1}A((TL')^{j+1})^{T}\right\}.
\end{align}
Then, it follows from lemma \ref{lem:TLeigen} that there exist $M>0$, $0<r<1$ such that
\begin{align}
&\sum_{j=0}^{t-2} \| (TL')^{j+1}A((TL')^{j+1})^{T} \|_{F} \leq \sum_{j=0}^{t-2} \| A\|_{F}\| (TL')^{j+1}\|^{2}_{F} \leq \sum_{j=0}^{t-2} {}^\exists M {}^\exists r^{j+1} \leq M \frac{r}{1-r}. \nonumber
\end{align}
Therefore, $S^{*}_{t}$ is absolutely convergent, so $S^{*}_{t}$ is convergent. We define this limit as $S^{*}_{\infty}:=\sum_{j=0}^{\infty} \left\{ (TL')^{j+1}A((TL')^{j+1})^{T}\right\}$ and prove that $S_{t}$ converges to $S^{*}_{\infty}$.

It follows that $\| S_{t}-S^{*}_{\infty} \|_{F} \leq \| S_{t} - S^{*}_{t}\|_{F} + \| S^{*}_{t} - S^{*}_{\infty} \|_{F}$, and $\| S^{*}_{t} - S^{*}_{\infty} \|_{F} \to 0$ because $S^{*}_{t} \rightarrow S^{*}_{\infty}$. Thus, if $\| S_{t} - S^{*}_{t}\|_{F} \to 0$, we get $\| S_{t}-S^{*}_{\infty} \|_{F} \to 0$. Now, $\| S_{t} - S^{*}_{t}\|_{F}$ can be bounded above as follows:
\begin{align}\label{eq:sa_St}
\| S_{t} - S^{*}_{t}\|_{F} &= \| \sum_{j=0}^{t-2} \left\{ (\prod_{i=0}^{j} TL_{t-i}')A_{t-j-1}(\prod_{i=0}^{j} TL_{t-i}')^{T} - (TL')^{j+1}A((TL')^{j+1})^{T} \right\} \|_{F} \nonumber
\\
&\leq  \sum_{j=0}^{t-2} \| \prod_{i=0}^{j} TL_{t-i}' - (TL')^{j+1}\|_{F} \| A_{t-j-1}\|_{F} \| (\prod_{i=0}^{j} TL_{t-i}')\|_{F} \nonumber
\\
&+ \sum_{j=0}^{t-2} \| (TL')^{j+1}\|_{F} \| A_{t-j-1}\|_{F} \| \prod_{i=0}^{j} TL_{t-i}' - (TL')^{j+1}\|_{F} \nonumber
\\
&+ \sum_{j=0}^{t-2} \| (TL')^{j+1}\|_{F} \| (A_{t-j-1}-A)\|_{F} \| (TL')^{j+1}\|_{F}.
\end{align}
For the first term of \eqref{eq:sa_St}, from the fact that $A_{t-j-1}$, $\prod_{i=0}^{j} TL_{t-i}'$ are bounded and lemma \ref{lem:sa_TL}, it follows that ${}^\exists M>0$, $0<{}^\exists r<1$ such that
\begin{align*}
&\sum_{j=0}^{t-2} \| \prod_{i=0}^{j} TL_{t-i}' - (TL')^{j+1}\|_{F} \| A_{t-j-1}\|_{F} \| (\prod_{i=0}^{j} TL_{t-i}')\|_{F} \leq r^{t}\sum_{j=0}^{t-2} (j+1)M \leq r^{t}t^{2}M.
\end{align*}
Then, $r^{t}t^{2}M \to 0 ~(t \to \infty)$, so the first term converges to 0. Similarly, the second term converges to 0. Finally, for the third term of \eqref{eq:sa_St}, using lemmas \ref{lem:TKFKT} and \ref{lem:TLeigen}, it follows that ${}^\exists M>0$, $0<{}^\exists r<1$ such that
\begin{align*}
&\sum_{j=0}^{t-2} \| (A_{t-j-1}-A)\|_{F} \| (TL')^{j+1}\|^{2}_{F} \leq \sum_{j=0}^{t-2} (M r^{t-j-1}) (M r^{j+1})^{2} \leq \sum_{j=0}^{t-2} M^{3} r^{t} = M^{3} (t-1) r^{t}.
\end{align*}
Then $M^{3} (t-1) r^{t} \to 0 ~(t \to \infty)$, so the third term converges to 0. Hence, $S_{t}$ converges to $S^{*}_{\infty}$, so $ZS_{t}Z^{T}$ that is the third term of \eqref{V_sum} converges to $ZS^{*}_{\infty}Z^{T}$. Thus, we get that all terms of \eqref{V_sum} converge. Therefore, $V(v_{t+1}')$ converges.
\end{proof}
%%%%%%%%%%%%%%%%%%%%%%%%%%%%%%%%%%%%%%%%%%%%%%%%%%%%%%%%%%%%%%%%%%%%
\section{Consistency}\label{sec:consistencyV}
In Section 3, $\hat V_{v'}=\frac{1}{n-n_{0}+1}\sum_{t=n_{0}}^{n} v_{t}'v_{t}'^T$ is an estimator of $V_{v'}$ that is a limiting value of $V(v'_{t})$. From Section 4, it follows that $E(v_{t}')\to 0$ and $V(v'_{t})\to V_{v'}$, so $E(\hat V_{v'})\to V_{v'}$. Thus, if $V\Bigl( \hat V_{v'}(i,j) \Bigl) =o(1)$, then $\hat V_{v'}$ is a consistent estimator of $V_{v'}$. We will show them in this section.

It follows that
\begin{align}\label{V_Vhat}
V\Bigl( \hat V_{v'}(i,j) \Bigl) &= V\Bigl( \frac{1}{n-n_{0}+1}\sum_{t=n_{0}}^{n} v_{i,t}'v_{j,t}' \Bigl) \nonumber
\\
&=\frac{1}{(n-n_{0}+1)^2}\Bigl\{ \sum_{t=n_{0}}^{n} V(v_{i,t}'v_{j,t}')+2\sum_{t=n_{0}}^{n-1}\sum_{s=1}^{n-t} Cov(v_{i,t}'v_{j,t}',\ v_{i,t+s}'v_{j,t+s}') \Bigl\} .
\end{align}
The first term of \eqref{V_Vhat} can be written as the following lemma:
\begin{lemma} \label{lem:Vv'v'}
$
\\
V(v_{i,t}'v_{j,t}') = V(v_{i,t}')V(v_{j,t}')+Cov(v_{i,t}',v_{j,t}')^{2}+E(v_{j,t}')^{2}V(v_{i,t}') + E(v_{i,t}')^{2}V(v_{j,t}') + 2E(v_{i,t}')E(v_{j,t}')Cov(v_{i,t}',\ v_{j,t}').
$
\end{lemma}
\noindent
\\
The proof is given in Appendix \ref{app:pf_lem}. Then $E(v_{i,t}')=o(1)$, $E(v_{j,t}')=o(1)$, $V(v_{i,t}')=O(1)$, $V(v_{j,t}')=O(1)$, and $Cov(v_{i,t}',v_{j,t}')=O(1)$ hold from the fact that $E(v_{t}')\to 0$ and $V(v'_{t})\to V_{v'}$. Thus, 
\begin{align*}
V(v_{i,t}'v_{j,t}')=&O(1).
\end{align*}
As a result, the first term of \eqref{V_Vhat} can be written as
\begin{align*}
\sum_{t=n_{0}}^{n} V(v_{i,t}'v_{j,t}')=O(n).
\end{align*}

Next, consider the second term of \eqref{V_Vhat}. We use the following lemma:
\begin{lemma} \label{lem:Covv'v'}
$
\\
Cov(v_{i,t}'v_{j,t}',\ v_{i,t+s}'v_{j,t+s}')=O(\|Cov(v_{t}',\ v_{t+s}')\|_{F}^{2} + \|Cov(v_{t}',\ v_{t+s}')\|_{F}).
$
\end{lemma}
The proof is given in Appendix \ref{app:pf_lem}. If $s=1$, then $Cov(v_{t}',\ v_{t+s}')$ can be written as
\begin{align*}
Cov(v_{t}',\ v_{t+1}')&=E \Bigl[ E \Bigl( (v_{t}' - E(v_{t}')) (v_{t+1}' - E(v_{t+1}'))^{T} |Y_{t} \Bigl) \Bigl]
\\
&=E \Bigl[ (v_{t}' - E(v_{t}')) E\Bigl( y_{t+1} - Za'_{t+1} - E(v_{t+1}') |Y_{t} \Bigl)^{T} \Bigl]
\\
&=E \Bigl[ (v_{t}' - E(v_{t}')) (Za_{t+1} - Za'_{t+1} - E(v_{t+1}'))^{T}  \Bigl] = Cov(v_{t}' ,\ a_{t+1} - a'_{t+1})Z^{T}.
\end{align*}
Thus, noting that the Frobenius norm is sub-multiplicative, we have
\begin{align*}
\| Cov(v_{t}',\ v_{t+1}')\|_{F} &\leq \| Cov(v_{t}' ,\ a_{t+1} - a'_{t+1})\|_{F} \| Z^{T}\|_{F}.
\end{align*}
Now we prepare the lemma as follows:
\begin{lemma} \label{lem:Cov_v}
$
\\
\text{Under assumptions 1,2,3,4, }{}^{\exists} M,\ {}^{\forall}t,\  \| Cov(v_{t}' ,\ a_{t+1} - a'_{t+1})\|_{F} <  M
$
\end{lemma}
\noindent
The proof is given in Appendix \ref{app:pf_lem}. Using lemma \ref{lem:Cov_v}, it follows that $\| Cov(v_{t}',\ v_{t+1}')\|_{F} \leq M \| Z^{T}\|_{F}$.
\\
\\
If $s\geq2$, then $Cov(v_{t}',\ v_{t+s}')$ can be written as
\begin{align*}
Cov(v_{t}',\ v_{t+s}')&=E \Bigl[  (v_{t}' - E(v_{t}')) E\Bigl( v_{t+s}' - E(v_{t+s}')  |Y_{t+s-1} \Bigl)^{T} \Bigl]
\\
&=E \Bigl[ (v_{t}' - E(v_{t}')) (Za_{t+s} - Za'_{t+s} - E(v_{t+s}'))^{T}  \Bigl]
\\
&=E \Bigl[ E \Bigl( (v_{t}' - E(v_{t}')) ( Za_{t+s} - Za'_{t+s} - E(v_{t+s}') )^{T}  | Y_{t+s-2} \Bigl) \Bigl]
\\
&=E \Bigl[ (v_{t}' - E(v_{t}')) E\Bigl( Za_{t+s} - Za'_{t+s} - E(v_{t+s}') | Y_{t+s-2} \Bigl)^{T} \Bigl].
\end{align*}
Now it holds (see the proof of lemma \ref{lem:E_a} in Appendix \ref{app:pf_lem}) that
\begin{align*}
E( a_{t+s} - a'_{t+s} | Y_{t+s-2} )=TL'_{t+s-1}(a_{t+s-1} - a'_{t+s-1}).
\end{align*}
Therefore,
\begin{align*}
&E \Bigl[ (v_{t}' - E(v_{t}')) E\Bigl( Za_{t+s} - Za'_{t+s} - E(v_{t+s}') | Y_{t+s-2} \Bigl)^{T} \Bigl]
\\
&=E\Bigl[ (v_{t}' - E(v_{t}')) ( ZTL'_{t+s-1}(a_{t+s-1} - a'_{t+s-1}) - E(v_{t+s}'))^{T} \Bigl] = Cov(v_{t}',\ a_{t+s-1} - a'_{t+s-1})(ZTL'_{t+s-1})^{T}.
\end{align*}
Using this transformation repeatedly, we have
\begin{align*}
Cov(v_{t}',\ v_{t+s}')&=Cov(v_{t}',\ a_{t+s-1} - a'_{t+s-1})(ZTL'_{t+s-1})^{T}
\\
&=\dots=Cov(v_{t}',\ a_{t+1} - a'_{t+1}) (\prod_{k=1}^{s-1} TL'_{t+s-k})^{T} Z^{T}.
\end{align*}
Additionally, by lemma \ref{lem:Cov_v} and lemma \ref{lem:sa_TL}, there exist $M>0$ and $0<r<1$ such that
\begin{align*}
\|Cov(v_{t}',\ v_{t+s}')\|_{F}&\leq \|  Cov(v_{t}',\ a_{t+1} - a'_{t+1})\|_{F} \| (\prod_{k=1}^{s-1} TL'_{t+s-k})^{T}\|_{F} \|Z^{T} \|_{F}\leq Msr^{s-1} \|Z^{T} \|_{F}.
\end{align*}
Note that this inequality holds if $s=1$. Thus, from the lemma \ref{lem:Covv'v'}, $Cov(v_{i,t}'v_{j,t}',\ v_{i,t+s}'v_{j,t+s}')$ can be bounded as follows:
\begin{align*}
Cov(v_{i,t}'v_{j,t}',\ v_{i,t+s}'v_{j,t+s}')&=O(\|Cov(v_{t}',\ v_{t+s}')\|_{F}^{2} + \|Cov(v_{t}',\ v_{t+s}')\|_{F})
\\
&=O( s^2r^{2(s-1)} M^{2} \|Z^{T}\|_{F}^{2} + sr^{s-1} M \|Z^{T}\|_{F})=(s^{2}r^{2(s-1)} + sr^{s-1})O(1).
\end{align*}
Then consider the second term of \eqref{V_Vhat}. Note that $0<r<1$; we have
\begin{align*}
\sum_{t=n_{0}}^{n-1}\sum_{s=1}^{n-t} Cov(v_{i,t}'v_{j,t}',\ v_{i,t+s}'v_{j,t+s}')&=\sum_{t=n_{0}}^{n-1}\sum_{s=1}^{n-t}(s^{2}r^{2(s-1)} + sr^{s-1})O(1)= \sum_{t=n_{0}}^{n-1}O(1)=O(n).
\end{align*}
From the above results, $V(\hat V_{v'}(i,j))$ can be bounded as follows:
\begin{eqnarray}
V(\hat V_{v'}(i,j))&=&\frac{1}{(n-n_{0}+1)^2}\{ \sum_{t=n_{0}}^{n} V(v_{i,t}'v_{j,t}')+2\sum_{t=n_{0}}^{n-1}\sum_{s=1}^{n-t} Cov(v_{i,t}'v_{j,t}',\ v_{i,t+s}'v_{j,t+s}')\} \nonumber
\\
&=&\frac{1}{(n-n_{0}+1)^2}(O(n)+O(n))=O(\frac{1}{n}).
\end{eqnarray}
Therefore, $V(\hat V_{v'}(i,j))\to 0$. Moreover, from $E(\hat V_{v'}(i,j))\to V_{v'}(i,j)$, which we mentioned at the beginning of this Section, it holds that
\begin{eqnarray}
\hat V_{v'}(i,j)\xrightarrow{p} V_{v'}(i,j).
\end{eqnarray}
That is, $\hat V_{v'}$ is a consistent estimator of $V_{v'}$.
%%%%%%%%%%%%%%%%%%%%%%%%%%%%%%%%%%%%%%%%%%%%%%%%%%%%%%%%%%%%%%%%%%%%
\section{Numerical example}\label{sec:sim1}
In Section 6, we conduct a Monte Carlo simulation and bus congestion forecasting to investigate the computational time and the forecast accuracy of our proposed algorithm. These simulations are run on Xeon Gold 6240R 2.4GHz with 512GB memory. We use R with the open BLAS library \citep{OpenBLAS} for fast matrix computation, and the package FKF \citep{FKF} to run the Kalman filter.
\subsection{Monte Carlo simulation}\label{sec:sim2}
\begin{figure}[H]
\centering
\includegraphics[width=12cm, bb=0 0 728 554]{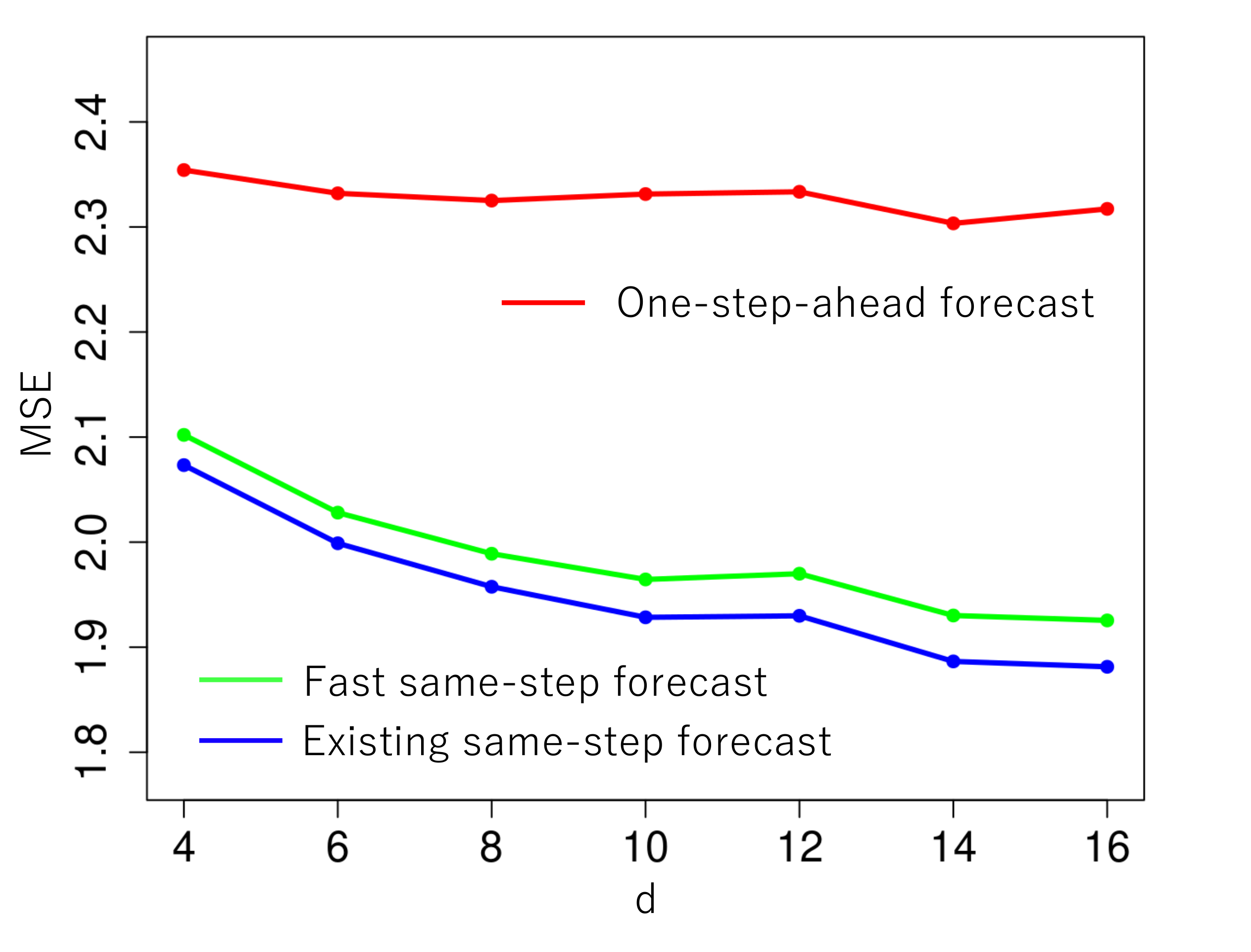}
\caption{Comparison of MSE between existing and fast method}
\label{fig:zu7}
\end{figure}

\begin{figure}[H]
\centering
\includegraphics[width=12cm, bb=0 0 738 558]{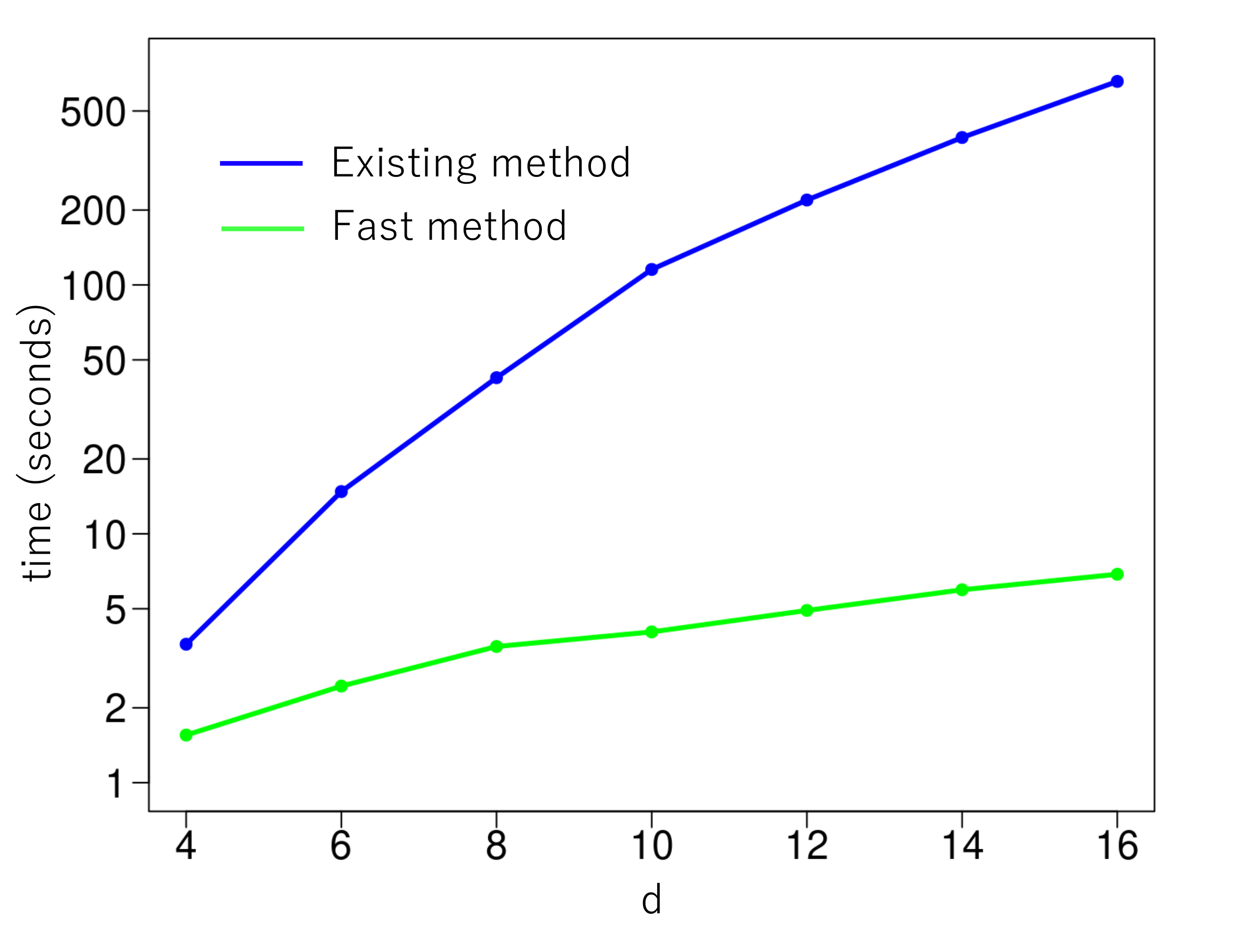}
\caption{Comparison of computation time between existing and fast method}
\label{fig:zu8}
\end{figure}

We compare the performance of the existing method in Section 2 and the fast method in Section 3 by Monte Carlo simulation. The simulation model is set to:
\begin{align} \label{eq:sim_mod}
\left\{ \begin{array}{ll}
 y_{j,t} = Z\alpha _{t}^{(j)}+\varepsilon _{j,t}   & \\
 \alpha _{t+1}^{(j)} = T\alpha _{t}^{(j)} + \eta _{t}^{(j)}  & \end{array} \right.  ( t = 1,2,\dots,2000, ~ j = 1,2,\dots, d), \nonumber
\end{align}
\begin{align*}
&Z=(1,\ 1,\ 0,\ 0,\ 0,\ 0,\ 0,\ 0),~
T=
\begin{pmatrix}
1 & 0 & 0 & 0 & 0 & 0 & 0 & 0\\
0 & -0.4 & -0.1 & 0 & 0 & 0 & 0.2 & 0.5 \\
0 & 1 & 0 & 0 & 0 & 0 & 0 & 0\\
0 & 0 & 1 & 0 & 0 & 0 & 0 & 0\\
0 & 0 & 0 & 1 & 0 & 0 & 0 & 0\\
0 & 0 & 0 & 0 & 1 & 0 & 0 & 0\\
0 & 0 & 0 & 0 & 0 & 1 & 0 & 0\\
0 & 0 & 0 & 0 & 0 & 0 & 1 & 0\\
\end{pmatrix}.
\end{align*}
Here, $\varepsilon _{t}\sim N\left( 0,\Sigma _{\varepsilon }\right)$, $\eta_{t}\sim N\left( 0,\Sigma_{\eta} \right)$, $\alpha _{1}^{(j)} = 0$, and
\begin{align*}
&\Sigma _{\varepsilon }(k,\ l) =
\left\{
\begin{array}{ll}
1 & (k=l)
\\
0.5 & (k\neq l),
\\
\end{array}
\right.
\\
&\Sigma _{\eta } = Diag(Q,\dots ,Q),~ Q(k,l)=
\left\{
\begin{array}{ll}
0.01 & (k=l=1)
\\
1 & (k=l=2)
\\
0 & (otherwise).
\\
\end{array}
\right.
\end{align*}
This is a combination of the AR(7) model and the local level model. In the simulation, we set the unknown parameters as $\theta_{1}:=\Sigma _{\varepsilon }(1,1),\dots , \theta_{d}:=\Sigma _{\varepsilon }(d,d)$, $\theta_{d+1}:=\Sigma _{\varepsilon }(k,l)~(k\neq l)$ and estimate these. The simulation also examines the effects of the changes in $d$ on the results. The following procedure is used in the simulation:
\begin{enumerate}
\item Generate the data $y_{1}, \dots ,y_{2000}$ that follow the above model and split these data as training data $Y_{train}:={y_{1},\dots , y_{1000}}$ and test data $Y_{test}:={y_{1001},\dots , y_{2000}}$.
\item Model estimation is performed using $Y_{train}$ by the existing and the fast methods, respectively. In the existing method, the maximum likelihood estimates for $\theta_{1},\dots ,\theta_{d+1}$ are computed. In the fast method, the maximum likelihood estimates for $\theta_{1},\dots ,\theta_{d}$ and also $\hat{V}_{v'}$ in \eqref{eq:vhat} are computed. Here, we select $n_{0} = 5$.
\item Run the same-step forecast of $y_{d, t}$ given $y_{1, t},\dots y_{d-1,t}$ for $t=1001,\dots ,2000$ by the existing method and the fast method, and then compute the squared forecast errors and measure the computation time taken to complete the estimation.
\end{enumerate}
For both existing and fast methods, the maximum likelihood estimates are computed numerically by the L-BFGS-B algorithm \citep{byrd1995limited}. We repeated this simulation 100 times, and then computed the mean-squared forecast error and the mean of computation time and compared them. Each dimension is $d=4,6,8,\dots ,16$. Figure \ref{fig:zu7} shows the mean-squared forecast error for the test data with 1000 observations over 100 runs. The results show that the same-step forecast of the fast method can modify the one-step-ahead forecast, similar to the existing method. The fast method is less accurate than the existing method, but not significantly less. Figure \ref{fig:zu8} shows the mean of computation time over 100 runs. The results show that the fast method is faster than the existing method for all $d=4,6,8,\dots ,16$. As $d$ is larger, the difference in computation time becomes larger; it could be because the computational complexity of the Kalman filter in the existing method is $O(d^{3})$, whereas in the fast method it is $O(d)$. In particular, when $d=16$, the fast method is approximately 90 times faster. The proposed method is slightly less accurate, but the computation time is much faster.

\subsection{Bus congestion forecasting}\label{sec:bus}
\begin{figure}[H]
\centering
\includegraphics[width=12cm, bb=0 0 737 594]{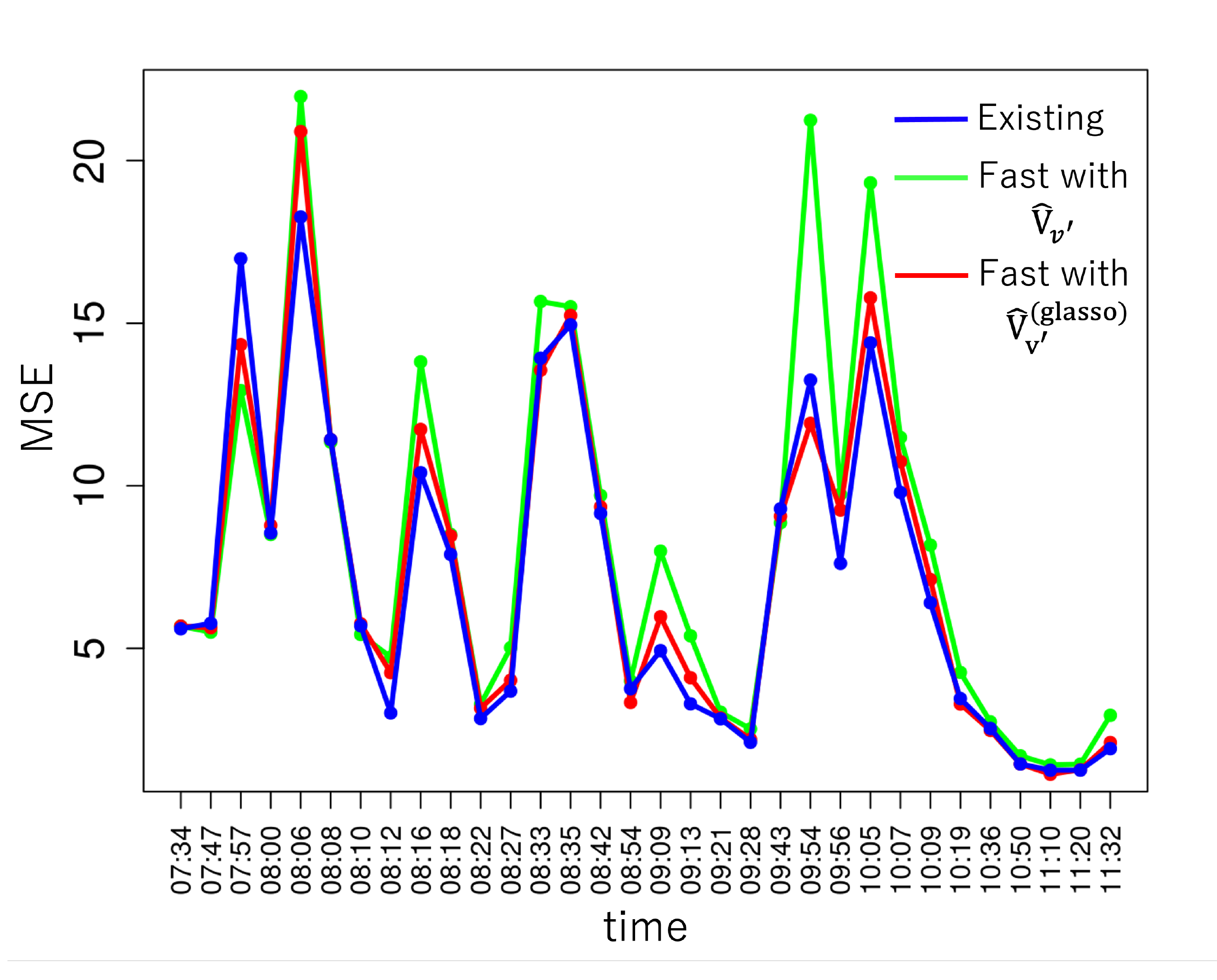}
\caption{Comparison of MSE of the same-step forecast between the existing and fast methods}
\label{fig:zu10}
\end{figure}

In Section \ref{sec:bus}, we investigate the effectiveness of the fast method through bus congestion forecasting. Such forecasts would be used as information to, for example, prevent COVID-19 infection. The data used are the numbers of passengers detected by sensors at the bus stops of ``Showa Bus'' \citep{sensys2020-takahashi}. This bus runs between Kyushu University and Kyudai-Gakkentoshi Station. Data from Monday to Friday mornings (32 buses per day), when the number of passengers is relatively large, are used. The data cover 379 days, with some missing values. The Kalman filter can interpolate missing values (e.g., see \citep{durbin2012time}). We let $y_{j,t}$ be the number of passengers detected by sensors of the $j$th bus on day $t$. Consider the following model:
\begin{align}
\left\{ \begin{array}{ll}
 y_{j,t} = Z\alpha _{t}^{(j)}+\varepsilon _{j,t}   & \\
 \alpha _{t+1}^{(j)} = T^{(j)}\alpha _{t}^{(j)} + \eta _{t}^{(j)} & 
 \end{array} \right.  ~~~( t = 1,2,\dots,379 ~~~ j = 1,2,\dots, 32), \nonumber
 \end{align}

\begin{align*}
&Z=(1,\ 1,\ 0,\ 0,\ 0,\ 0),~
T^{(j)}=
\begin{pmatrix}
1 & 0 & 0 & 0 & 0 & 0 \\
0 & \phi_{1}^{(j)} & \phi_{2}^{(j)} & \phi_{3}^{(j)} & \phi_{4}^{(j)} & \phi_{5}^{(j)} \\
0 & 1 & 0 & 0 & 0 & 0 \\
0 & 0 & 1 & 0 & 0 & 0 \\
0 & 0 & 0 & 1 & 0 & 0 \\
0 & 0 & 0 & 0 & 1 & 0 \\
\end{pmatrix}.
\end{align*}
Here, $\varepsilon _{t}\sim N\left( 0,\Sigma _{\varepsilon }\right)$, $\eta_{t}\sim N\left( 0,\Sigma_{\eta} \right)$, $\alpha _{1}\sim N( 0 , 10^{7} \times I_{192})$, and
\begin{align*}
&\Sigma _{\eta } = Diag(Q^{(1)},\dots ,Q^{(32)}),~  Q^{(j)}(k,l)=
\left\{
\begin{array}{ll}
q_{1}^{(j)} & (k=l=1)
\\
q_{2}^{(j)} & (k=l=2)
\\
0 & (otherwise).
\\
\end{array}
\right.
\end{align*}
This is a combination of the AR(5) model and the local level model, as in the simulation in Section \ref{sec:sim2}. The difference is that the initial value of the state vector $\alpha _{1}$ is unknown. Therefore, we assume that $\alpha _{1}$ follows uninformative prior $N( 0 , 10^{7} \times I_{192})$. The unknown parameters are $\phi_{1}^{(j)},\dots , \phi_{5}^{(j)}$, $q_{1}^{(j)}$, $q_{2}^{(j)}$ $(j = 1,\dots ,32)$, and all components of $\Sigma _{\varepsilon }$. The number of parameters is $7\times 32 + \frac{32\times 33}{2} = 752$. The existing method requires the estimation of all these parameters. In contrast, the fast method does not require estimation of the off-diagonal components of $\Sigma _{\varepsilon }$. Therefore, the number of parameters to be estimated is $8\times 32 = 256$. The following procedure is used in the simulation:
\begin{enumerate}
\item Split the data as training data $Y_{train}:={y_{1},\dots , y_{190}}$ and test data $Y_{test}:={y_{191},\dots , y_{379}}$.
\item  Model estimation is performed using $Y_{train}$ by the existing and fast methods. In the existing method, the maximum likelihood estimates for $\phi_{1}^{(j)},\dots , \phi_{5}^{(j)}$, $q_{1}^{(j)}$, $q_{2}^{(j)}$ $(j = 1,\dots ,32)$ and $\Sigma _{\varepsilon }$ are computed. In the fast method, the maximum likelihood estimates for $\phi_{1}^{(j)},\dots , \phi_{5}^{(j)}$, $q_{1}^{(j)}$, $q_{2}^{(j)}$, $\Sigma _{\varepsilon }(j,j)$ $(j = 1,\dots ,32)$ and $\hat{V}_{v'}$ or $\hat V_{v'}^{(glasso)}$ are computed. Here, we select $n_{0} = 5$, and we use parallel computing for the estimation for each $j = 1,\dots, 32$. The number of parallels is 32; thus, they are all computed in parallel. The regularization parameter $\lambda$ in $\hat V_{v'}^{(glasso)}$ is determined by the BIC \citep{schwarz1978estimating}.
\item Run the same-step forecast of $y_{j, t}$ for $t=191,\dots ,379$, $j = 1,\dots ,32$ by the existing method and the fast method, and then compute the mean-squared forecast error and measure the computation time taken to complete the estimation. Here, the same-step forecast uses data up to 10 minutes before; we may use data 1--9 minutes before to achieve high accuracy, but the passengers may not need such a very short-term forecast value.
\end{enumerate}
For both existing and fast methods, the maximum likelihood estimates are computed numerically by the L-BFGS-B algorithm. The log transformation of $Y_{train}$ is performed to stabilize model estimation, but the original $Y_{test}$ is used to calculate the prediction error. Figure \ref{fig:zu10} shows the mean-squared forecast error of the same-step forecast at each time. The fast method with $\hat{V}_{v'}$ is slightly less accurate than other methods. The fast method with $\hat V_{v'}^{(glasso)}$ is comparable to the existing method. The computation time is $4.001 \times 10^{4}$ seconds (about 11 hours) for the existing method, 6.951 seconds for the fast method with $\hat{V}_{v'}$, and 7.905 seconds for the fast method with $\hat V_{v'}^{(glasso)}$, respectively.

From the results of this data analysis, when $d$ is large and many unknown parameters have to be estimated, as in this bus congestion data, the fast method is significantly faster and thus practical. Also, when $n$ is not large enough compared with $d$, using $\hat V_{v'}^{(glasso)}$ for the fast same-step forecast is better, and its forecast accuracy would be comparable to the existing method.

%%%%%%%%%%%%%%%%%%%%%%%%%%%%%%%%%%%%%%%%%%%%%%%%%%%%%%%%%%%%%%%%%%%%
\section{Concluding remarks}
In this paper, we considered the same-step forecast and proposed the fast method. In this method, we estimated the mean vector and the covariance matrix of the one-step-ahead forecast error. The estimation is based on their convergence, and we discussed the estimators' consistency. A Monte Carlo simulation was conducted to investigate the effectiveness of the fast method. Bus congestion forecasting was also performed to illustrate the usefulness of the fast method. The results showed that our proposed method was much faster than the existing method. In bus congestion forecasting, our method was comparable to the existing method.

In future studies, we will consider researching the convergence of the covariance matrix of the one-step-ahead forecast error when $Z_{t}$ or $T_{t}$ is time-varying. This convergence can ensure that when $n$ is large, the same-step forecast provides better forecast accuracy than the one-step-ahead forecast in the fast method. Therefore, the discussion of the convergence would make it more clear in what cases of $Z_{t}$ and $T_{t}$ the fast same-step forecast improves the fast one-step-ahead forecast. Finally, we will also consider researching a method to select the best $n_{0}$ in \eqref{eq:vhat} that minimizes the mean-squared same-step forecast error. From our observations, we believe it is beneficial to set a small value for $n_{0}$, but we are not yet confident which value of $n_{0}$ is best.
%%%%%%%%%%%%%%%%%%%%%%%%%%%%%%%%%%%%%%%%%%%%%%%%%%%%%%%%%%%%%%%%%%%%
%%%%%%%%%%%%%%%%%%%%%%%%%%%%%%%%%%%%%%%%%%%%%%%%%%%%%%%%%%%%%%%%%%%%
\section*{Appendix}
\addcontentsline{toc}{section}{Appendix}

\appendix

\section{Proof of lemmas}\label{app:pf_lem}

\begin{flalign*}
&\text{\bf{Lemma \ref{lem:E_a} } }&
\\
&E(a_{t+1}-a_{t+1}') =  TL_{t}'E(a_{t}-a_{t}') &
\end{flalign*}

\begin{proof}
Generally, we consider time-varying $Z_{t}$, $T_{t}$ and prove $E(a_{t+1}-a_{t+1}') =  T_{t}L_{t}'E(a_{t}-a_{t}')$ instead of this lemma.
\\
Using \eqref{eq:kalman}, \eqref{eq:kalman_nocor2}, $E(v_{t}|Y_{t-1})=0$, and the fact that $a_{t}$, $a_{t}'$ are constant when $y_{1}, . . . , y_{t-1}$ are given, we have
\begin{align*}
&E(a_{t+1}-a_{t+1}')=E(T_{t}a_{t}+T_{t}K_{t}v_{t}-T_{t}a_{t}'-T_{t}K_{t}'v_{t}')
\\
&=T_{t} E \Bigl( E(a_{t}+K_{t}v_{t}-a_{t}'-K_{t}'v_{t}'|Y_{t-1}) \Bigl) = T_{t} E \Bigl( E(a_{t}-a_{t}'-K_{t}'(y_{t}-Z_{t}a_{t}')|Y_{t-1})  \Bigl)
\\
&=T_{t} E \Bigl( a_{t}-a_{t}'-K_{t}'(Z_{t}a_{t}-Z_{t}a_{t}') \Bigl)=T_{t}(I_{p}-K_{t}'Z_{t})E(a_{t}-a_{t}')=T_{t}L_{t}'E(a_{t}-a_{t}').
\end{align*}
\end{proof}
\begin{flalign*}
&\text{\bf{Lemma \ref{lem:V_a} } }&
\\
&E \Bigl( (a_{t+1}-a_{t+1}')(a_{t+1}-a_{t+1}')^{T} \Bigl) = T(K_{t}-K_{t}')F_{t}(K_{t}-K_{t}')^{T}T^{T}+TL_{t}'E \Bigl( (a_{t}-a_{t}')(a_{t}-a_{t}')^{T} \Bigl) L_{t}'^{T}T^{T}&
\end{flalign*}

\begin{proof}
Generally, we consider time-varying $Z_{t}$, $T_{t}$ and prove $E \Bigl( (a_{t+1}-a_{t+1}')(a_{t+1}-a_{t+1}')^{T} \Bigl) = T_{t}(K_{t}-K_{t}')F_{t}(K_{t}-K_{t}')^{T}T_{t}^{T}+T_{t}L_{t}'E \Bigl( (a_{t}-a_{t}')(a_{t}-a_{t}')^{T} \Bigl) L_{t}'^{T}T_{t}^{T}$ instead of this lemma.
\\
Note $v_{t}'=y_{t}-Z_{t}a_{t}'=y_{t}-Z_{t}a_{t}+Z_{t}a_{t}-Z_{t}a_{t}'=v_{t}+Z_{t}a_{t}-Z_{t}a_{t}'$, we have
\begin{align*}
&E \Bigl( (a_{t+1}-a_{t+1}')(a_{t+1}-a_{t+1}')^{T} \Bigl)
\\
&=E \Bigl( (T_{t}a_{t}+T_{t}K_{t}v_{t}-T_{t}a_{t}'-T_{t}K_{t}'v_{t}')(T_{t}a_{t}+T_{t}K_{t}v_{t}-T_{t}a_{t}'-T_{t}K_{t}'v_{t}')^{T} \Bigl)
\\
&=T_{t} E \Bigl( (a_{t}+K_{t}v_{t}-a_{t}'-K_{t}'(v_{t}+Z_{t}a_{t}-Z_{t}a_{t}')) (a_{t}+K_{t}v_{t}-a_{t}'-K_{t}'(v_{t}+Z_{t}a_{t}-Z_{t}a_{t}'))^{T} \Bigl) T_{t}^{T}
\\
&=T_{t} E \Bigl( ( (I_{p}-K_{t}'Z_{t})(a_{t}-a_{t}')+(K_{t}-K_{t}')v_{t})( (I_{p}-K_{t}'Z_{t})(a_{t}-a_{t}')+ (K_{t}-K_{t}')v_{t})^{T} \Bigl) T_{t}^{T}
\\
&=T_{t} \Bigl\{ E \Bigl( L_{t}'(a_{t}-a_{t}')(a_{t}-a_{t}')^{T}L_{t}'^{T} \Bigl) + E\Bigl( L_{t}'(a_{t}-a_{t}')v_{t}^{T}(K_{t}-K_{t}')^{T} \Bigl)
\\
& ~~~~~~~~~~~~~~~~~~~~~~~~~~~+ E\Bigl( (K_{t}-K_{t}')v_{t}(a_{t}-a_{t}')^{T}L_{t}'^{T} \Bigl) + E\Bigl( (K_{t}-K_{t}')v_{t}v_{t}^{T}(K_{t}-K_{t}')^{T} \Bigl) \Bigl\} T_{t}^{T}.
\end{align*}
Now, considering the second term, since $a_{t}$, $a_{t}'$ are constant when $y_{1}, . . . , y_{t-1}$ are given, it follows that
\begin{align*}
& E\Bigl( L_{t}'(a_{t}-a_{t}')v_{t}^{T}(K_{t}-K_{t}')^{T} \Bigl)
\\
&=E \Bigl[ E\Bigl( L_{t}'(a_{t}-a_{t}')v_{t}^{T}(K_{t}-K_{t}')^{T}|Y_{t-1} \Bigl) \Bigl] =E \Bigl[ L_{t}'(a_{t}-a_{t}')E(v_{t}^{T}|Y_{t-1})(K_{t}-K_{t}')^{T} \Bigl]
\\
&=E \Bigl[ L_{t}'(a_{t}-a_{t}')\times 0\times (K_{t}-K_{t}')^{T} \Bigl] =0.
\end{align*}
Similarly, it follows that $E\Bigl( (K_{t}-K_{t}')v_{t}(a_{t}-a_{t}')^{T}L_{t}'^{T} \Bigl) =0$. Considering the fourth term, we have
\begin{align*}
&E\Bigl( (K_{t}-K_{t}')v_{t}v_{t}^{T}(K_{t}-K_{t}')^{T} \Bigl)
\\
&=(K_{t}-K_{t}')E(v_{t}v_{t}^{T})(K_{t}-K_{t}')^{T}=(K_{t}-K_{t}')E\Bigl( E(v_{t}v_{t}^{T}|Y_{t-1}) \Bigl) (K_{t}-K_{t}')^{T}.
\end{align*}
Since $E(v_{t}|Y_{t-1})=0$, it follows that $E\Bigl( E(v_{t}v_{t}^{T}|Y_{t-1}) \Bigl)=E\Bigl( V(v_{t}|Y_{t-1}) \Bigl)=E(F_{t})=F_{t}$. Hence, it holds that
\begin{align*}
&E\Bigl( (K_{t}-K_{t}')v_{t}v_{t}^{T}(K_{t}-K_{t}')^{T} \Bigl) =(K_{t}-K_{t}')F_{t}(K_{t}-K_{t}')^{T}.
\end{align*}
From the above, we get
\begin{align*}
&E \Bigl( (a_{t+1}-a_{t+1}')(a_{t+1}-a_{t+1}')^{T} \Bigl)
\\
&=T_{t}(K_{t}-K_{t}')F_{t}(K_{t}-K_{t}')^{T}T_{t}^{T}+T_{t}L_{t}'E \Bigl( (a_{t}-a_{t}')(a_{t}-a_{t}')^{T} \Bigl) L_{t}'^{T}T_{t}^{T}.
\end{align*}
\end{proof}

\begin{flalign*}
&\text{\bf{Lemma \ref{lem:sa_TL}} }&
\\
&\text{Under assumptions 1,2,3,4, }{}^\exists M>0,\ 0<{}^\exists r<1,\ \| \prod_{i=0}^{j} TL'_{t-i} - (TL')^{j+1} \|_{F}\leq (j+1)Mr^{t}, &
\\
&{}^\exists M >0,\ 0<{}^\exists r<1,\ \| \prod_{i=0}^{j} TL'_{t-i} \|_{F} \leq M(j+2)r^{j+1}. &
\end{flalign*}

\begin{proof}
For $j \geq 1$, it follows that
\begin{align*}
&\| \prod_{i=0}^{j} TL'_{t-i} - (TL')^{j+1} \|_{F} \leq \| (TL'_{t} - TL') \prod_{i=1}^{j} TL'_{t-i}\|_{F} + \| TL' (\prod_{i=1}^{j} TL'_{t-i} - (TL')^{j}) \|_{F}
\\
&\leq \dots
\\
&\leq \sum_{k=0}^{j-1} \left\{ \| (TL')^{k} (TL'_{t-k} - TL') \prod_{i=k+1}^{j} TL'_{t-i}\|_{F}\right\} + \| (TL')^{j} (TL'_{t-j} - TL')\|_{F}
\\
&\leq \sum_{k=0}^{j-1} \left\{ \| (TL')^{k}\|_{F} \| (TL'_{t-k} - TL')\|_{F} \| \prod_{i=k+1}^{j} TL'_{t-i}\|_{F}\right\} + \| (TL')^{j}\|_{F} \| (TL'_{t-j} - TL')\|_{F}
\\
&\leq \sum_{k=0}^{j-1} \left [ {}^\exists M_{1} {}^\exists r^{k}_{1}\right ] \left [ \| T\|_{F} {}^\exists M_{2} {}^\exists r^{t-k}_{2} \right ] {}^\exists M_{3} + M_{1} r^{j}_{1}  M_{2}  r^{t-j}_{2} ~~~~~(\because \text{lemma } \ref{lem:TLeigen}, \text{ lemma } \ref{lem:PFKLcon},\text{ and assumption }4)
\\
&\leq \sum_{k=0}^{j-1}  (M r^{t})  + M r^{t} = (j+1)Mr^{t},
\end{align*}
where $M:=M_{1}M_{2} \max{ \{M_{3} \| T\|_{F} , 1\} }$, $r:=\max{ \{r_{1}, r_{2} \} }$. Note that this inequality holds for $j=0$ from lemma $\ref{lem:PFKLcon}$. In addition, note $j+1 \leq t$, by this inequality and lemma \ref{lem:TLeigen}, it holds that there exist $M>0$ and $0<r<1$, such that
\begin{align*}
\| \prod_{i=0}^{j} TL'_{t-i}\|_{F} &\leq \| \prod_{i=0}^{j} TL'_{t-i} - (TL')^{j+1} \|_{F} + \|(TL')^{j+1}\|_{F}
\\
&\leq (j+1)Mr^{t} + Mr^{j+1} \leq (j+1)Mr^{j+1} + Mr^{j+1} = M(j+2)r^{j+1}.
\end{align*}
\end{proof}

\begin{flalign*}
&\text{\bf{Lemma \ref{lem:con_TL}} }&
\\
&\text{Under assumptions 1,2,3,4, }\prod_{i=0}^{t-1} TL_{t-i}'\to 0(t\to \infty ). &
\end{flalign*}
\begin{proof}
From lemma \ref{lem:sa_TL}, it holds that there exist $M>0$ and $0<r<1$, such that
\begin{align*}
\| \prod_{i=0}^{t-1} TL'_{t-i}\|_{F} \leq M(t+1)r^{t} \rightarrow 0(t\rightarrow \infty).
\end{align*}
\end{proof}

\begin{flalign*}
&\text{\bf{Lemma \ref{lem:TKFKT}} }&
\\
& \text{Under assumptions 1,2,3, }{}^\exists M>0,\ 0<{}^\exists r<1,\ &
\\
& \| T(K_{t}-K_{t}')F_{t}(K_{t}-K_{t}')^{T}T^{T} - T(K-K')F(K-K')^{T}T^{T}\|_{F} \leq Mr^{t}. &
\end{flalign*}

\begin{proof}
Note $F_{t}$ and $K_{t}-K_{t}'$ are finite, we have
\begin{align*}
&\| T(K_{t}-K_{t}')F_{t}(K_{t}-K_{t}')^{T}T^{T} - T(K-K')F(K-K')^{T}T^{T}\|_{F} 
\\
&\leq \| T\|^{2}_{F} \Big{\{} \|(K_{t}-K_{t}') - (K-K')\|_{F} \|F_{t}\|_{F} \| (K_{t}-K_{t}') \|_{F} 
\\
&+ \|(K-K')\|_{F} \|F_{t} - F\|_{F} \| (K_{t}-K_{t}') \|_{F} + \|(K-K')\|_{F} \|F\|_{F} \| (K_{t}-K_{t}') - (K-K')\|_{F} \Big{\}}
\\
&\leq {}^\exists M_{1}(2\|(K_{t}-K_{t}') - (K-K')\|_{F} + \|F_{t} - F\|_{F})
\\
&\leq {}^\exists M_{1}(2\|K_{t}-K\|_{F} + 2\|K'_{t}-K'\|_{F} + \|F_{t} - F\|_{F}) \leq {}^\exists M {}^\exists r^{t} ~~~~~~(\because \text{lemma } \ref{lem:PFKLcon}).
\end{align*}
\end{proof}

\begin{flalign*}
&\text{\bf{Lemma \ref{lem:Vv'v'}} }&
\\
&V(v_{i,t}'v_{j,t}') = V(v_{i,t}')V(v_{j,t}')+Cov(v_{i,t}',v_{j,t}')^{2}+E(v_{j,t}')^{2}V(v_{i,t}') + E(v_{i,t}')^{2}V(v_{j,t}') + 2E(v_{i,t}')E(v_{j,t}')Cov(v_{i,t}',\ v_{j,t}'). &
\end{flalign*}

\begin{proof}
$V(v_{i,t}'v_{j,t}')$ can be written as follows:
\begin{align*}
V(v_{i,t}'v_{j,t}')=&V \Bigl( (v_{i,t}'-E(v_{i,t}'))(v_{j,t}'-E(v_{j,t}')) + (v_{i,t}'-E(v_{i,t}'))E(v_{j,t}') + E(v_{i,t}')(v_{j,t}'-E(v_{j,t}')) + E(v_{i,t}')E(v_{j,t}') \Bigl)
\\
=&V \Bigl( (v_{i,t}'-E(v_{i,t}'))(v_{j,t}'-E(v_{j,t}')) \Bigl) + V \Bigl( (v_{i,t}'-E(v_{i,t}'))E(v_{j,t}') + E(v_{i,t}')(v_{j,t}'-E(v_{j,t}')) \Bigl)
\\
&+ 2Cov \Bigl( (v_{i,t}'-E(v_{i,t}'))(v_{j,t}'-E(v_{j,t}')),\ (v_{i,t}'-E(v_{i,t}'))E(v_{j,t}') + E(v_{i,t}')(v_{j,t}'-E(v_{j,t}')) \Bigl).
\end{align*}
Consider the covariance of the third term in the last equation. Since the mean of the right side is 0, we have
\begin{align*}
&Cov \Bigl( (v_{i,t}'-E(v_{i,t}'))(v_{j,t}'-E(v_{j,t}')),\ (v_{i,t}'-E(v_{i,t}'))E(v_{j,t}') + E(v_{i,t}')(v_{j,t}'-E(v_{j,t}')) \Bigl)
\\
&=E \Bigl( ((v_{i,t}'-E(v_{i,t}'))(v_{j,t}'-E(v_{j,t}')) \Bigl[ (v_{i,t}'-E(v_{i,t}'))E(v_{j,t}') + E(v_{i,t}')(v_{j,t}'-E(v_{j,t}')) \Bigl] \Bigl) =0.
\end{align*}
In the last equality, we use the fact that a third order moment of a normal random variable with $0$ mean is $0$. Therefore,
\begin{align*}
V(v_{i,t}'v_{j,t}')=&V \Bigl( (v_{i,t}'-E(v_{i,t}'))(v_{j,t}'-E(v_{j,t}')) \Bigl) + V \Bigl( (v_{i,t}'-E(v_{i,t}'))E(v_{j,t}') + E(v_{i,t}')(v_{j,t}'-E(v_{j,t}')) \Bigl)
\\
= &E \Bigl( \{ (v_{i,t}'-E(v_{i,t}'))(v_{j,t}'-E(v_{j,t}'))\} ^{2} \Bigl) - E \Bigl( (v_{i,t}'-E(v_{i,t}'))(v_{j,t}'-E(v_{j,t}')) \Bigl)^{2}
\\
&+ E(v_{j,t}')^{2}V(v_{i,t}') + E(v_{i,t}')^{2}V(v_{j,t}') + 2E(v_{i,t}')E(v_{j,t}')Cov(v_{i,t}',\ v_{j,t}')
\\
= &E \Bigl( \{ (v_{i,t}'-E(v_{i,t}'))(v_{j,t}'-E(v_{j,t}'))\} ^{2} \Bigl) - Cov(v_{i,t}',v_{j,t}')^{2}
\\
&+ E(v_{j,t}')^{2}V(v_{i,t}') + E(v_{i,t}')^{2}V(v_{j,t}') + 2E(v_{i,t}')E(v_{j,t}')Cov(v_{i,t}',\ v_{j,t}').
\end{align*}
Moreover, using Isserlis' theorem \citep{isserlis1918formula}, we have
\begin{align*}
&E \Bigl( \{ (v_{i,t}'-E(v_{i,t}'))(v_{j,t}'-E(v_{j,t}'))\} ^{2} \Bigl)
\\
&=E \Bigl( (v_{i,t}'-E(v_{i,t}'))(v_{j,t}'-E(v_{j,t}')) (v_{i,t}'-E(v_{i,t}'))(v_{j,t}'-E(v_{j,t}')) \Bigl) =V(v_{i,t}')V(v_{j,t}')+2Cov(v_{i,t}',v_{j,t}')^{2}.
\end{align*}
Therefore, we get
\begin{align*}
&V(v_{i,t}'v_{j,t}') = V(v_{i,t}')V(v_{j,t}')+Cov(v_{i,t}',v_{j,t}')^{2}+E(v_{j,t}')^{2}V(v_{i,t}') + E(v_{i,t}')^{2}V(v_{j,t}') + 2E(v_{i,t}')E(v_{j,t}')Cov(v_{i,t}',\ v_{j,t}').
\end{align*}
\end{proof}

\begin{flalign*}
&\text{\bf{Lemma \ref{lem:Covv'v'}} }&
\\
&Cov(v_{i,t}'v_{j,t}',\ v_{i,t+s}'v_{j,t+s}')=O(\|Cov(v_{t}',\ v_{t+s}')\|_{F}^{2} + \|Cov(v_{t}',\ v_{t+s}')\|_{F}).&
\end{flalign*}

\begin{proof}
Similar to the proof of lemma \ref{lem:Vv'v'}, we have
\begin{align*}
&Cov(v_{i,t}'v_{j,t}',\ v_{i,t+s}'v_{j,t+s}')
\\
=&Cov \Bigl( (v_{i,t}'-E(v_{i,t}'))(v_{j,t}'-E(v_{j,t}')),\  (v_{i,t+s}'-E(v_{i,t+s}'))(v_{j,t+s}'-E(v_{j,t+s}')) \Bigl)
\\
&+Cov \Bigl( (v_{i,t}'-E(v_{i,t}'))E(v_{j,t}') + E(v_{i,t}')(v_{j,t}'-E(v_{j,t}')),
\\
&~~~~~~~~~~~~~~~~~~~~~~~~~~~~~~~~~~~~~~(v_{i,t+s}'-E(v_{i,t+s}'))E(v_{j,t+s}') + E(v_{i,t+s}')(v_{j,t+s}'-E(v_{j,t+s}')) \Bigl)
\\
=&E \Bigl( (v_{i,t}'-E(v_{i,t}'))(v_{j,t}'-E(v_{j,t}'))(v_{i,t+s}'-E(v_{i,t+s}'))(v_{j,t+s}'-E(v_{j,t+s}')) \Bigl)
\\
&-E \Bigl( (v_{i,t}'-E(v_{i,t}'))(v_{j,t}'-E(v_{j,t}')) \Bigl) E \Bigl( (v_{i,t+s}'-E(v_{i,t+s}'))(v_{j,t+s}'-E(v_{j,t+s}')) \Bigl)
\\
&+E(v_{j,t}')E(v_{j,t+s}')Cov(v_{i,t}',\ v_{i,t+s}')+E(v_{i,t}')E(v_{i,t+s}')Cov(v_{j,t}',\ v_{j,t+s}')
\\
&+E(v_{j,t}')E(v_{i,t+s}')Cov(v_{i,t}',\ v_{j,t+s}')+E(v_{i,t}')E(v_{j,t+s}')Cov(v_{j,t}',\ v_{i,t+s}')
\\
=&Cov(v_{i,t}',\ v_{i,t+s}') Cov(v_{j,t}',\ v_{j,t+s}') + Cov(v_{i,t}',\ v_{j,t+s}') Cov(v_{j,t}',\ v_{i,t+s}') + Cov(v_{i,t}',\ v_{j,t}') Cov(v_{i,t+s}',\ v_{j,t+s}')
\\
& - Cov(v_{i,t}',\ v_{j,t}') Cov(v_{i,t+s}',\ v_{j,t+s}')
\\
&+E(v_{j,t}')E(v_{j,t+s}')Cov(v_{i,t}',\ v_{i,t+s}')+E(v_{i,t}')E(v_{i,t+s}')Cov(v_{j,t}',\ v_{j,t+s}')
\\
&+E(v_{j,t}')E(v_{i,t+s}')Cov(v_{i,t}',\ v_{j,t+s}')+E(v_{i,t}')E(v_{j,t+s}')Cov(v_{j,t}',\ v_{i,t+s}')~~~~(\because \text{Isserlis' theorem).}
\\
=&Cov(v_{i,t}',\ v_{i,t+s}') Cov(v_{j,t}',\ v_{j,t+s}') + Cov(v_{i,t}',\ v_{j,t+s}') Cov(v_{j,t}',\ v_{i,t+s}')
\\
&+E(v_{j,t}')E(v_{j,t+s}')Cov(v_{i,t}',\ v_{i,t+s}')+E(v_{i,t}')E(v_{i,t+s}')Cov(v_{j,t}',\ v_{j,t+s}')
\\
&+E(v_{j,t}')E(v_{i,t+s}')Cov(v_{i,t}',\ v_{j,t+s}')+E(v_{i,t}')E(v_{j,t+s}')Cov(v_{j,t}',\ v_{i,t+s}')
\\
\leq&\ 2 \|Cov(v_{t}',\ v_{t+s}')\|_{F}^{2} 
\\
&+ \Bigl\{ |E(v_{j,t}')E(v_{j,t+s}')|+|E(v_{i,t}')E(v_{i,t+s}')|+|E(v_{j,t}')E(v_{i,t+s}')|+|E(v_{i,t}')E(v_{j,t+s}')| \Bigl\} \|Cov(v_{t}',\ v_{t+s}')\|_{F}.
\end{align*}
From this inequality, it follows that
\begin{align*}
&Cov(v_{i,t}'v_{j,t}',\ v_{i,t+s}'v_{j,t+s}')=O(\|Cov(v_{t}',\ v_{t+s}')\|_{F}^{2} + \|Cov(v_{t}',\ v_{t+s}')\|_{F}).
\end{align*}
\end{proof}

\begin{flalign*}
&\text{\bf{Lemma \ref{lem:Cov_v}} }&
\\
&\text{Under assumptions 1,2,3,4, }{}^{\exists} M,\ {}^{\forall}t,\  \| Cov(v_{t}' ,\ a_{t+1} - a'_{t+1})\|_{F} <  M &
\end{flalign*}

\begin{proof}
It follows that
\begin{align*}
&\| Cov(v_{t}' ,\ a_{t+1} - a'_{t+1})\|_{F} = \| Cov(y_{t} - Za_{t}' ,\ Ta_{t}+TK_{t}v_{t} - Ta_{t}'-TK_{t}'v_{t}')\|_{F} 
\\
&= \| Cov(v_{t} + Z(a_{t} - a_{t}') ,\  T(K_{t}-K_{t}')v_{t} + T(I_{p}-K_{t}'Z)(a_{t}-a_{t}') )\|_{F}
\\
&= \| Cov(v_{t} + Z(a_{t} - a_{t}') ,\  T(K_{t}-K_{t}')v_{t} + TL_{t}'(a_{t}-a_{t}') )\|_{F}.
\end{align*}
From the fact that $a_{t}$, $a_{t}'$ are constant when $y_{1}, . . . , y_{t-1}$ are given, it follows that
\begin{align*} 
Cov(v_{t},\ a_{t} - a_{t}') = E\Bigl( v_{t} (a_{t} - a_{t}' - E(a_{t} - a_{t}') )^{T} \Bigl) &=E\Bigl( E( v_{t} | y_{1},\dots ,y_{t-1}) (a_{t} - a_{t}' - E(a_{t} - a_{t}') )^{T} \Bigl)
\\
&=E\Bigl( 0 \times (a_{t} - a_{t}' - E(a_{t} - a_{t}') )^{T} \Bigl) =0.
\end{align*}
Thus, note that $V(v_{t}) = F_{t}$; we have
\begin{align*}
&\| Cov(v_{t} + Z(a_{t} - a_{t}') ,\  T(K_{t}-K_{t}')v_{t} + TL_{t}'(a_{t}-a_{t}') )\|_{F}
\\
&=\| F_{t} (T(K_{t}-K_{t}'))^{T} + Z V(a_{t} - a_{t}') (TL_{t}')^{T}\|_{F}
\\
&\leq \| F_{t} (T(K_{t}-K_{t}'))^{T}\|_{F} + \|Z\|_{F} \| V(a_{t} - a_{t}')\|_{F} \| (TL_{t}')^{T}\|_{F}
\\
&\leq \| F_{t} (T(K_{t}-K_{t}'))^{T}\|_{F} + \|Z\|_{F} \Bigl\{ \| E\Bigl( (a_{t} - a_{t}')(a_{t} - a_{t}')^{T}\Bigl) \|_{F} + \| E(a_{t} - a_{t}')E(a_{t} - a_{t}')^T \|_{F} \Bigl\} \| (TL_{t}')^{T}\|_{F}.
\end{align*}
Now, since $F_{t}$, $K_{t}$, and $K_{t}'$ converge, $F_{t} (T(K_{t}-K_{t}'))^{T}$ also converge. Therefore, it follows that $^{\exists}M_{1},\  ^{\forall}t,\ \| F_{t} (T(K_{t}-K_{t}'))^{T}\|_{F}<M_{1}$. Similarly, since  $L_{t}'$ converge, it follows that $^{\exists}M_{2},\  ^{\forall}t,\ \| (TL_{t}')^{T}\|_{F}<M_{2}$, and since $E\Bigl( (a_{t} - a_{t}')(a_{t} - a_{t}')^{T}\Bigl) $ and $E(a_{t} - a_{t}')$ converge (see Section \ref{sec:convergenceEV}), it follows that ${}^{\exists}M_{3},\ {}^{\exists}M_{4},\  ^{\forall}t,\ \| E\Bigl( (a_{t} - a_{t}')(a_{t} - a_{t}')^{T}\Bigl) \|_{F}<M_{3},\ \| E(a_{t} - a_{t}')E(a_{t} - a_{t}')^T \|_{F}<M_{4}$. From the above, let $M := M_{1} + \|Z\|_{F} \times M_{2}(M_{3}+M_{4})$, for any $t$ it holds that
\begin{align*}
&\| Cov(v_{t}' ,\ a_{t+1} - a'_{t+1})\|_{F} = \| Cov(v_{t} + Z(a_{t} - a_{t}') ,\  T(K_{t}-K_{t}')v_{t} + TL_{t}'(a_{t}-a_{t}') )\|_{F} < M.
\end{align*}
\end{proof}

\section{Supplement of Lemma \ref{lem:TLeigen}}\label{app:TL'geocon}
We prove that for any matrix A of which all eigenvalues are of absolute value less than 1 it holds that $^\exists M>0$, $0<^\exists r<1$, $\| A^{n}\|_{F} \leq Mr^{n}$. First, it follows that $A^{n} \rightarrow 0$ from that all eigenvalues of the matrix A are of absolute value less than 1. Then, there exists a positive constant c large enough, such that
\begin{align*}
&\| A^{c} \|_{F} < {}^\exists r_{1} < 1.
\end{align*}
Let q be the remainder when n is divided by c, for all $n \in \mathbb{N}$ it follows that
\begin{align*}
&\| A^{n} \|_{F} = \| A^{ c \times \lfloor \frac{n}{c} \rfloor + q} \|_{F} \leq r_{1}^{\lfloor \frac{n}{c} \rfloor} \| A^{q} \|_{F} 
\end{align*}
Now $r_{1}^{\lfloor \frac{n}{c} \rfloor}$ can be bounded above as
\begin{align*}
& r_{1}^{\lfloor \frac{n}{c} \rfloor} = \frac{1}{r_{1}} r_{1}^{1 + \lfloor \frac{n}{c} \rfloor} \leq \frac{1}{r_{1}} r_{1}^{ \frac{n}{c} },
\end{align*}
and let $M_{1} := \max_{ q } \| A^{q} \|_{F} $, we get 
\begin{align*}
&\| A^{n} \|_{F} \leq M_{1} \frac{1}{r_{1}} r_{1}^{\frac{n}{c} } = Mr^{n},
\end{align*}
where $M:=M_{1} \frac{1}{r_{1}}$, $r:=r_{1}^{ \frac{1}{c} }$. Note $r<1$ from $r_{1}<1$, the inequality we want was given.

\section{Mild assumption for the convergence of $E(v_{t}')$}\label{app:conE_assum}
In Section 4, we proved the convergence of $E(v_{t}')$ and $V(v_{t}')$ under assumption \ref{assumption} and the condition that $Z_{t}=Z$, $T_{t}=T$. However, the convergence of $E(v_{t}')$ even holds under the condition that $Z_{t}=Z$, $T_{t}=T$ are time-varying, and the assumptions can be mild as follows:
\begin{align*}
&\text{Assumption for the convergence of $E(v_{t}')$ : }\prod_{i=0}^{t-1} T_{t-i}L_{t-i}'\to 0~~~(t\to \infty ).
\end{align*}
Note this assumption is more mild than assumption \ref{assumption} from lemma \ref{lem:con_TL}. We prove the convergence of $E(v_{t}')$ under this assumption.
\begin{proof}
Since lemma \ref{lem:E_a} holds even if $Z_{t}=Z$, $T_{t}=T$ are time-varying from Appendix \ref{app:pf_lem}, it follows that
\begin{align*}
&E(v_{t+1}')=E(y_{t+1}-Z_{t}a_{t+1}')=E\Bigl( E(y_{t+1}-Z_{t}a_{t+1}'|Y_{t}) \Bigl) =E(Z_{t}a_{t+1}-Z_{t}a_{t+1}')
\\
&=Z_{t}E(a_{t+1}-a_{t+1}')=Z_{t}T_{t}L_{t}'E(a_{t}-a_{t}')=\dots =Z_{t}(\prod_{i=0}^{t-1} T_{t-i}L_{t-i}')(a_{1}-a_{1}')~~~(\because \text{lemma \ref{lem:E_a}}).
\end{align*}
Since $\prod_{i=0}^{t-1} T_{t-i}L_{t-i}'\to 0(t\to \infty )$ from the assumption above, it holds that $Z_{t}(\prod_{i=0}^{t-1} T_{t-i}L_{t-i}')(a_{1}-a_{1}')\to 0(t\to \infty )$.
Therefore, we get $E(v_{t}')\to 0(t\to \infty )$. 
\end{proof}
%%%%%%%%%%%%%%%%%%%%%%%%%%%%%%%%%%%%%%%%%%%%%%%%%%%%%%%%%%%%%%%%%%%%%%%%%%%%%%%%%%%%%%%%%%%%%%%%%%%%%%%%%%%%%%%%%%%%%%%%%%%%%%%%%%%%%%%%%%%%%%%%%%%%%%%%%%%%%%%%%%%%%%%%%%%%%%%%%%%%%%%%%%%%%%%%%%%%%%%%%%

\section*{Acknowledgments}
The authors would like to thank Arakawa Lab for providing us with the bus sensor data. This work was supported by JSPS KAKENHI Grant Number 22J20435 and WISE program (MEXT) at Kyushu University.

\renewcommand{\refname}{References}

\bibliographystyle{plainnat}
\bibliography{reference_combine}

\end{document}